\documentclass[final]{svjour3}
\usepackage{euscript,amsmath,bm,amssymb,amsfonts,graphicx,epstopdf}

\usepackage{color}
\usepackage{soul}
\usepackage{mathtools}
\numberwithin{equation}{section}

\newcommand{\E}{{\mathbb E}}

\renewcommand{\u}{\widetilde{u}}

\newcommand{\e}{{\rm e}}
\renewcommand{\r}{{\bf r}}

\newcommand{\hxi}{{\bm \xi}}

\newcommand{\w}{\bm{w}}

\newcommand{\x}{\bm{x}}

\newcommand{\y}{\bm{y}}

\newcommand{\calU}{\mathcal{U}}

\newcommand{\calA}{\mathcal{A}}

\newcommand{\calN}{\mathcal{N}}

\newcommand{\calJ}{\mathcal{J}}

\newcommand{\calI}{\mathcal{I}}

\newcommand{\calG}{\mathcal{G}}

\newcommand{\R}{\mathbb R}

\newcommand{\n}{\bm{n}}

    \numberwithin{equation}{section}

\newcommand{\abs}[1]{\vert #1 \vert}

\newcommand{\Markov}[2]{\underset{#1}{\overset{#2}{\rightleftharpoons}}}

\begin{document}

\title{\bf Cellular diffusion processes in singularly perturbed domains}

\author{Paul C. Bressloff}

\institute{Paul C. Bressloff \at Department of
  Mathematics, Imperial College London, \\ London SW7 2AZ UK
  \\
  \email{p.bressloff@imperial.ac.uk} }

\authorrunning{Bressloff}

\maketitle

    \begin{abstract}
There are many processes in cell biology that can be modeled in terms of particles diffusing in a two-dimensional (2D) or three-dimensional (3D) bounded domain $\Omega \subset \R^d$ containing a set of small subdomains or interior compartments $\calU_j$, $j=1,\ldots,N$ (singularly-perturbed diffusion problems). The domain $\Omega$ could represent the cell membrane, the cell cytoplasm, the cell nucleus or the extracellular volume, while an individual compartment could represent a synapse, a membrane protein cluster, a biological condensate, or a quorum sensing bacterial cell. In this review we use a combination of matched asymptotic analysis and Green's function methods to solve a general type of singular boundary value problems (BVP) in 2D and 3D, in which an inhomogeneous Robin condition is imposed on each interior boundary $\partial \calU_j$. This allows us to incorporate a variety of previous studies of singularly perturbed diffusion problems into a single mathematical  modeling framework. We mainly focus on steady-state solutions and the approach to steady-state, but also highlight some of the current challenges in dealing with time-dependent solutions and randomly switching processes. \end{abstract}

    \section{Introduction}
    
 There are a wide range of problems in cell biology that can be modelled in terms of particles diffusing in a two-dimensional (2D) or three-dimensional (3D) bounded domain $\Omega \subset \R^d$ containing a set of small subdomains or compartments $\calU_j$, $j=1,\ldots,N$, within the interior (singularly-perturbed diffusion problems). For example, $\Omega$ could represent the cell membrane while the compartments could represent synapses or protein clusters such as macro-assemblies of adhesion molecules. Alternatively, if $\Omega$ were the cell volume then the corresponding compartments could be enzyme-rich reactive substrates, or high density liquid droplets (biological condensates). On a larger length scale, $\Omega$ could represent the extracellular medium surrounding a set of diffusively-coupled cells as in bacterial quorum sensing. In most of these applications one is interested in calculating the steady-state solution (assuming it exists), and the rate of approach to steady state as determined by the principal eigenvalue of the linear evolution operator or by the so-called accumulation time\footnote{Another class of problem treats the compartments as totally or partially absorbing traps, and the main focus is determining the first passage time or splitting probability for a single particle to be captured by one of the traps -- the so-called narrow capture problem \cite{Coombs09,Chevalier11,Coombs15,Ward15,Lindsay16,Lindsay17,Grebenkov20,Bressloff21a,Bressloff21b}. We do not consider single particle diffusion in this review.}.
 The latter is based on the fractional deviation of the time-dependent solution from steady state, and has been used extensively in studies of diffusion-based morphogenesis \cite{Berez10,Berez11,Berez11a,Gordon11,Bressloff19}.
The type of steady-state boundary value problem (BVP) that typically arises can be analysed using a combination of matched asymptotic analysis and Green's function methods. The basic idea is to construct an inner or local solution of the diffusion equation  in a small neighbourhood of each compartment, which is then matched to an outer solution in the bulk domain. The matching is achieved by expressing the outer solution in terms of the Green's function of the diffusion equation in the absence of any compartments. It follows that the details of the matched asymptotic analysis in 2D and 3D domains differ considerably due to corresponding differences in the Green's function singularities: 
\begin{align*}
&G(\x,\x_0)\rightarrow -\frac{1}{2\pi D}\ln|\x-\x_0| \mbox{ in 2D },  \\ &G(\x,\x_0)\rightarrow \frac{1}{4\pi D|\x-\x_0|} \mbox{ in 3D}
\end{align*}
as $|\x-\x_0|\rightarrow 0$.
Hence, an asymptotic expansion of the solution to a BVP in 3D is in powers of $\epsilon$, where $\epsilon$ represents the size of a compartment relative to the size of the bulk domain. On the other hand, the analogous expansion in 2D tends to be in powers of $\nu=-1/\ln \epsilon$ at $O(1)$ in $\epsilon$. The slower convergence of $\nu$ in the limit $\epsilon \rightarrow 0$ can be handled by summing the logarithmic terms non-perturbatively, as originally shown in Refs. \cite{Ward93,Ward93a}. 

Examples of cellular processes that have been analyzed using singular perturbation theory include the following: intracellular protein concentration gradients mediated by enzyme-rich subcellular domains \cite{Straube09,Levy11,Levy13}; coarsening of biological condensates via Ostwald ripening \cite{Kavanagh14,Bressloff20a,Bressloff20b,Bressloff24b}; protein clustering \cite{Bressloff24b}; synaptic receptor trafficking \cite{Bressloff08,Bressloff23a}; bacterial quorum sensing \cite{Muller13,Gou16,Iyan21a,Iyan21b}; volume transmission \cite{Lawley18,Lawley20}. Although each of these applications utilises the same basic mathematical method, they tend to be treated in isolation. Therefore, in this review we show how the majority of current applications can be incorporated into a single type of BVP involving an inhomogeneous Robin boundary condition on $\partial \calU_j$ combined with  the insertion and removal of particles within $\Omega$. We consider three different models for the boundary inhomogeneity $c_j$: (I) $c_j$ is a prescribed constant on $\partial \calU_j$; (II) $\partial \calU_j$ acts as a semi-permeable membrane and $c_j$ is identified with the particle concentration on the interior side of the membrane; (iii)  $\calU_j$ is treated as a well-mixed biochemical compartment that can exchange one molecular species with the exterior environment., whose concentration within $\calU_j$ determines $c_j$. For simplicity, we assume that all of the compartments are either circularly symmetric (2D) or spherically symmetric (3D). More general shapes are briefly discussed in section 6.

The structure of the paper is as follows\footnote{In this review we focus on the underlying mathematical principles and techniques, and how they are applied to singularly-perturbed diffusion problems in cell biology. Explicit results for specific configurations of $\Omega$ and $\calU_j$, $j=1,\ldots,N$, and choices of boundary conditions can be found in the cited literature.}. In section 2 we introduce the general model equations and show how they are related to the various applications in specific limits. In section 3 we solve the steady-state BVP in both 2D and 3D using singular perturbation theory. The relaxation to steady-state is analysed in section 4 by solving the BVP for the singularly perturbed diffusion equations in Laplace space and then taking the small-$s$ limit, where $s$ is the Laplace variable. This determines the corresponding accumulation time. In section 5 we summarise a recent study \cite{Iyan21b} that uses singular perturbation theory to reduce a PDE-ODE model of diffusively coupled biochemical compartments (model III)  to a nonlinear system of ODEs. Finally, in section 6, we describe various extensions of the theory, including randomly switching boundaries, and Kuramoto-like models of diffusively coupled oscillators. These examples highlight one major challenge for future work on singularly-perturbed diffusion problems, namely, developing asymptotic and numerical methods that can deal effectively with time-dependent and switching processes. (For very recent progress in this direction see Ref. \cite{Pelz24}.)
 
     \setcounter{equation}{0}

\section{The general setup}

\begin{figure*}[t!]
\centering
\includegraphics[width=12cm]{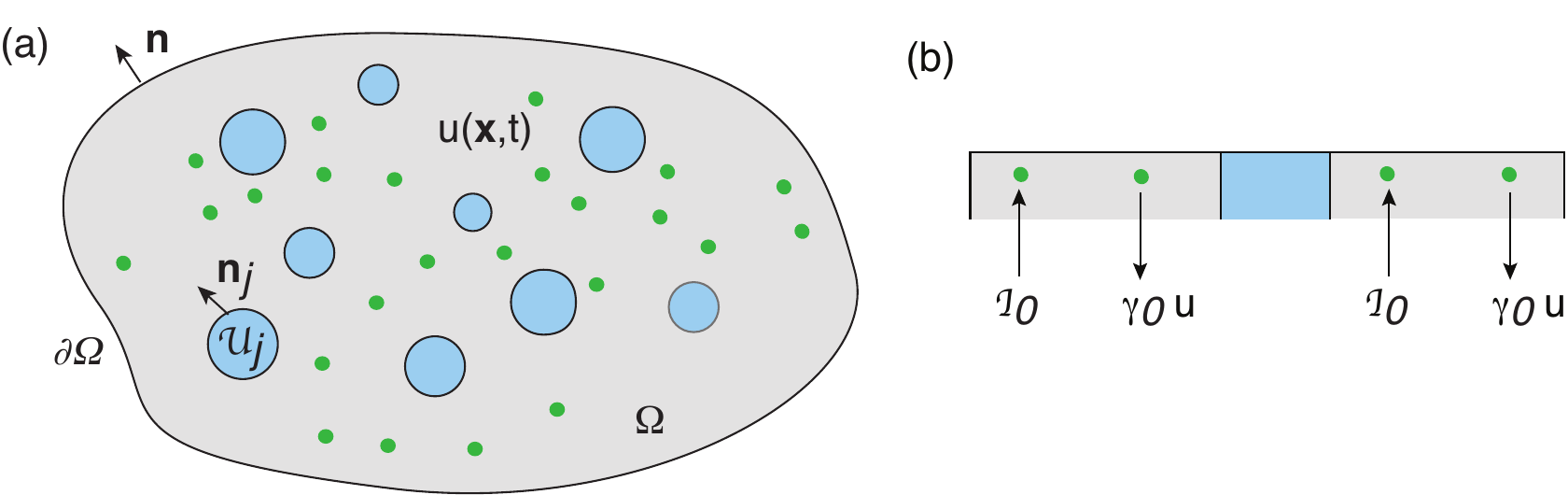} 
\caption{(a) A singularly perturbed domain $\Omega\subset \R^d$ containing $N$ interior compartments $\calU_j$ centred at the positions $\x_j$, $j=1,\ldots,N$. (b) In the region exterior to the compartments, $\Omega\backslash \calU$ with $\calU=\bigcup_{j=1}^N\calU_j$, there may be the insertion and/or removal of diffusing particles. The influx current is denoted by $\calI_0$, whereas the efflux current is taken to be proportional to the local particle concentration $u(\x,t)$ with rate $\gamma_0$.}
\label{fig1}
\end{figure*}

Consider a $d$-dimensional simply-connected, bounded domain $\Omega \subset \R^d$, $d=2,3$, containing a set of $N$ interior subdomains or compartments $\calU_j$, $j=1,\ldots,N$, see Fig. \ref{fig1}(a). For concreteness, suppose that each compartment is circularly (spherically) symmetric for $d=2$ ($d=3$) -- more general shapes will be considered in section 6. Denoting the radius and centre of the $j$-th compartment by $r_j$ and $\x_j$, respectively, we have
\begin{equation}
\calU_j=\{\x\in \Omega,\ |\x-\x_j|< r_j\},\quad \partial \calU_j=\{\x\in \Omega,\ |\x-\x_j|= r_j\}.
\end{equation}
The main characteristic of a singularly perturbed domain is that the compartments $\calU_j$ are small compared to the size of the domain $\Omega$ and are well separated. More precisely, suppose that the domain $\Omega$ is inscribed by a rectangular area or volume whose smallest dimension is $L$, and introduce the dimensionless parameter
$\epsilon =r_{\max}/L$ where $r_{\max}=\max \{r_j,\, j=1,\ldots,N\}$. We then fix length scales by setting $L=1$ and writing $r_j=\epsilon \ell_j$ with $\ell_j=r_j/r_{\max}$,
and $0<\epsilon \ll 1$. We also assume that $|\x_i-\x_j| =O(1)$ for all $j\neq i$ and $\min_{{\bf s}}\{|\x_j -{\bf s}|,{\bf s} \in \partial \Omega \} =O(1)$, $j=1,\ldots,N$. (If $\Omega$ is unbounded, then we identify $L$ with the smallest distance between any pair of compartments, that is, $L=\min\{|\x_i-\x_j|,\ i\neq j\}$.)

\begin{figure}[t!]
\centering
\includegraphics[width=12cm]{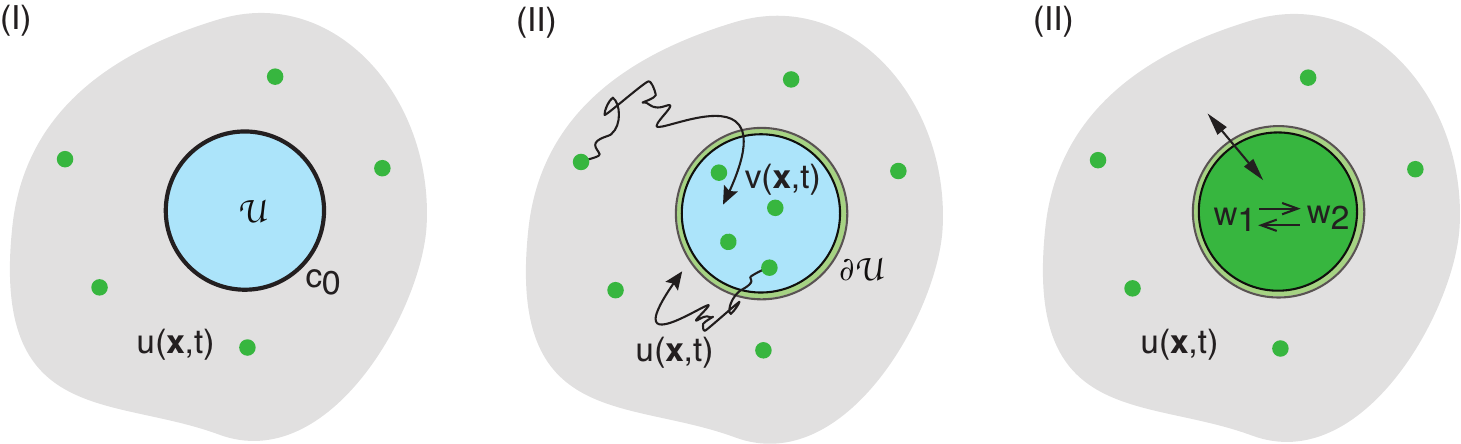} 
\caption{Three distinct models of the boundary field $c(\x,t)$, $\x \in \partial \calU$, of a compartment $\calU$, see equation (\ref{gmaster}c). (I) The boundary field $c(\x,t)$ for $\x\in\partial \calU$ is a fixed constant $c_0$. (II) The boundary $\partial \calU$ acts as a semipermeable membrane and particles can freely diffuse across the interface. The boundary field $c(\x,t)=v(\x^-,t)$ for $\x\in \calU$, where $\x^-$ indicates that the point on the boundary is approached from the interior. (III) $\calU$ is treated as a well-mixed biochemical compartment containing $K$ chemical species $w_a$, $a=1,\ldots,K$. The boundary field $c(\x,t)=w_1(t)$ for all $\x\in \partial \calU$.}
\label{fig2}
\end{figure}

Let $u(\x,t)$ denote the particle concentration in the region exterior to all of the compartments, that is, $\x \in \Omega\backslash \calU$ with $\calU=\bigcup_{j=1}^N\calU_j$. We take $u$ to satisfy the following singularly perturbed diffusion problem:
\begin{subequations} 
\label{gmaster}
\begin{align}
	\frac{\partial u(\x,t)}{\partial t} &= D{\bm \nabla}^2u(\x,t)-\gamma_0 u(\x,t) +\calI_0,\quad \x\in \Omega \backslash   \calU,\\
	D{\bm \nabla} u(\x,t) \cdot \n(\x)&=0,\  \x\in\partial \Omega ,\\
	D  {\bm \nabla} u(\x,t)\cdot \n_j(\x) &=\kappa_j[u(\x,t) -c_j(\x,t)],\quad \x\in \partial \calU_j,	\end{align}
	\end{subequations}
	where $D$ is the diffusivity, $\n$ is the outward unit normal of the exterior surface $\partial \Omega$ and $\n_j$ is the unit normal of the interior surface $\partial \calU_j$ that is directed out of the interior $\calU_j$, see Fig. \ref{fig1}(a).
The term $\calI_0$ in equation (\ref{gmaster}a) represents a particle source current, whereas the term $-\gamma_0u(\x,t)$ represents removal of particles at a rate $\gamma_0$, see Fig. \ref{fig1}(b); the latter could correspond to a degradation rate. The exterior boundary $\partial \Omega$ is assumed to be totally reflecting, whereas on each interior boundary $\partial \calU_j$ we impose a generalised Robin boundary condition with constant reactivity $\kappa_j$. In the limit $\kappa_j\rightarrow \infty$, equation (\ref{gmaster}c) reduces to the inhomogeneous Dirichlet boundary condition $u(\x,t)=c_j(\x,t)$ for all $\x\in \partial \calU_j$, whereas we obtain a totally reflecting boundary condition when $\kappa_j=0$. It remains to specify the boundary field $c_j(\x,t)$, $\x\in \partial \calU_j$. We consider three different compartmental models for $c_j(\x,t)$, as illustrated in  Fig. \ref{fig2}:
\medskip

\noindent {\em I. Constant boundary fields}. The simplest boundary condition is $c_j(\x,t)=c_{j,0}$ for all $\x\in \partial \calU_j$, where $c_{j,0}$ is a prescribed constant. Equation (\ref{gmaster}c) then reduces to a classical Robin boundary condition for $u(\x,t)-c_{j,0}$. \medskip

\noindent {\it II. Spatially inhomogeneous compartments and semi-permeable interfaces.}  Now suppose that the boundary $\partial \calU_j$ acts as a semipermeable interface between diffusion in the exterior domain $\Omega\backslash\calU$ and diffusion within each interior $\calU_j$. Equations (\ref{gmaster}a,b) are supplemented by the following diffusion equation for the particle concentration $v_j(\x,t)$ in $\calU_j$:
\begin{subequations}
\label{gmasterj}
\begin{align}
	\frac{\partial v_j(\x,t)}{\partial t}&= \overline{D}_j {\bm \nabla}^2 v_j(\x,t)  -\overline{\gamma}_j v_j(\x,t)+\overline{\calI}_j, \ \x\in \calU_j, 
\end{align}
and the boundary condition (\ref{gmaster}c) is replaced by the pair of semipermeable boundary conditions
\begin{align}
D{\bm \nabla} u(\x^+,t)\cdot \n_j(\x^+) &= \overline{D}_j{\bm \nabla} v_j(\x^-,t)\cdot \n_j(\x^-) =\kappa_j[u(\x^+,t) -v_j(\x^-,t)] ,\ \x  \in \partial \calU_j.
\end{align}
	\end{subequations}
Here $\x^{\pm}$ indicates whether a point on the boundary $\partial \calU_j$ is approached from the exterior or interior of the compartment $\calU_j$.
The diffusivity $\overline{D}_j$ within the $j$-th compartment may differ from the exterior diffusivity $D$. We also allow for the insertion and removal of particles within $\calU_j$. \medskip 

\noindent {\em III. Spatially homogeneous compartments as biochemical reaction networks.}  Each compartment is treated as a biochemical reaction network involving a set of $K+1$ chemical species $X_a$, $a=0,1,\ldots,K$. Diffusion within $\calU_j$ is assumed to be sufficiently fast so that the corresponding compartment can be treated as well-mixed. This means that the concentrations $w_{j,a}$ within $\calU_j$ are spatially homogeneous. Only one of the chemical species can be exchanged with the exterior domain $\Omega\backslash \calU$, which we identify with $X_0$. This implies that $c_j(\x,t)\equiv w_{j,0}(t)$ for all $\x\in \partial \calU_j$. Using mass action kinetics we then have the following system of equations for the concentrations ${\bf w}_j=(w_{j,0},w_{j,0},\ldots,w_{j,K})$:
\begin{align}
\label{mass}
\frac{dw_{j,a}}{dt}=f_a(\w_j)-  \frac{\kappa_j \delta_{a,1}}{|\calU_j|} \int_{\partial \calU_j}( w_{j,0}(t)- u(\x,t))d\x
\end{align}
for $ a=1,\dots,K $. For simplicity, the mass action kinetics is taken to be  the same for all cells, that is, the form of $f_a$ is independent of $j$.) One can view equations (\ref{gmaster}) and (\ref{mass}) with $c_j(\x,t)=w_{j,0}(t)$ as a nonlinear PDE-ODE model of diffusion-mediated communication between small signalling compartments \cite{Gou16,Iyan21a,Iyan21b}.
Finally, note that the factor of $|\calU_j|$ ensures that the number of particles crossing $\partial \calU_j$ is consistent with the definition of the flux in the boundary condition (\ref{gmaster}c), assuming that the compartmental and bulk concentrations are defined with respect to the same volume.

We now summarise various biological applications of equations (\ref{gmaster}) that involve one of the above compartmental models. 
\medskip

\begin{figure}[t!]
  \centering
  \includegraphics[width=5cm]{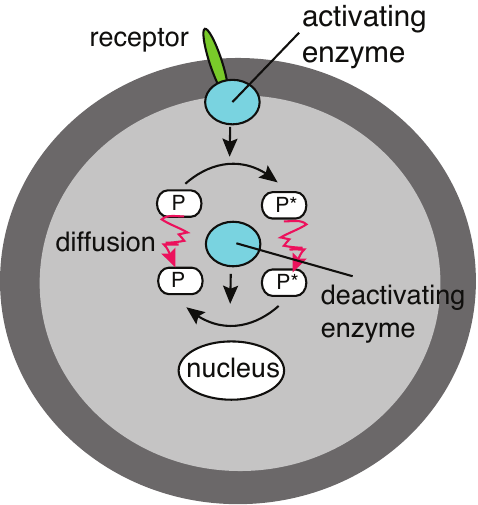}
  \caption{Schematic illustration of a protein modification cycle where an inactive form $P$ is converted to an active form $P^*$  at the plasma membrane. Both forms diffuse in the cytoplasm, resulting in deactivation of $P^*$ by cytoplasmic enzymes.}
  \label{protgrad}
\end{figure}

\noindent \underline{Intracellular signalling \cite{Straube09,Levy11,Levy13}.} An important component of many intracellular signal transduction pathways is the reversible cycling between an inactive and an active protein state, which is catalysed by opposing activator and deactivator enzymes. A concentration gradient in these signalling cycles can then be generated by the spatial segregation of the opposing enzymes \cite{Brown99,Kholodenko09,Garcia10,Howard12}. One such mechanism is the phosphorylation of proteins by a membrane-bound kinase, which are then dephosphorylated by a cytosolic phosphatase, see Fig. \ref{protgrad}. This results in a gradient of the phosphorylated protein, with a high concentration close to the cell membrane and a low concentration within the interior of the cell. As the cell grows in size, the surface-to-volume ratio decreases and membrane-activated proteins have to diffuse over longer distances in order to reach their target such as the nucleus. Hence, the proteins become progressively deactivated towards the cell interior, thus providing a mechanism for coupling cell growth with the cell cycle and cell division \cite{Meyers06}. It has also been observed that the range of an intracellular protein gradient can be enhanced by signalling cascades mediated by kinases located on the membranes of subcellular compartments. Cascades of intracellular phosphorylation-dephosphorylation gradients can be formulated in terms of equations (\ref{gmaster}) for $\calI_0=0$ and $c_j(\x,t)=c_{j,0}$ (model I) \cite{Straube09,Levy11,Levy13}. The boundaries $\partial \calU_j$ represent the membranes of subcellular compartments that are rich in kinases, and $u(\x,t)$ is the concentration of an activated signalling molecule in the cytosol. Phosphotases within the cytosol deactivate the signalling molecules at the constant rate $ \gamma_0$. If the total cytosolic concentration $u_{\rm tot}$ of activated and deactivated signalling molecules is a constant and the enzymes on $\partial \calU_j$ are unsaturated, then one can identify $c_{j,0}$ with $u_{\rm tot}$ and $\kappa_j$ with the rate of phosphorylation. 

\medskip

\noindent \underline{Ostwald ripening of biological condensates \cite{Alikakos04,Kavanagh14,Bressloff20a,Bressloff20b,Bressloff24b}.} Ostwald ripening describes the coarsening of polydisperse droplets during late-stage liquid-liquid phase separation, ultimately transitioning to a single condensate in thermodynamic equilibrium with a surrounding dilute phase \cite{Lif61,Wagner61}. In this example $\calU_j$ represents a liquid droplet in the high density phase, $u(\x,t)$ is the concentration of particles in the low density phase, and $\gamma_0=0=\calI_0=0$. The boundary condition on each surface is taken to be $u(\x,t)=c_{j,0}$ with
\begin{align}
\label{GT0}
c_{j,0}=\phi_a \left (1+\frac{\ell_c}{r_j}\right )\equiv c_0(r_j),
\end{align}
where $\phi_a$ is the density of the dilute phase under complete phase separation, and $\ell_c$ is the so-called capillary length constant. (More precisely, $\phi_a$ denotes a volume fraction.)  
Equation (\ref{GT0}) expresses the Gibbs-Thomson law due to interfacial tension on the droplet interface \cite{Doi13}. We thus have an example of model I in the limit $\kappa_j\rightarrow \infty$. The no-flux boundary condition on $\partial \Omega$ ensures mass conservation. Suppose that the coarsening dynamics is much slower than the equilibration of the concentration profile in the dilute phase. Under this quasi-static approximation, the solute concentration $u(\x)$ exterior to the droplets satisfies a steady-state diffusion equation of the form
\begin{equation}
 \label{Mullins}
{\bm \nabla}^2 u=0,\quad \x \in \Omega\backslash \calU,\quad {\bm \nabla} u(\x)\cdot \n(\x) =0 \mbox{ on } \partial \Omega,\quad u(\x)=c_{0}(r_j)  \mbox{ on }\partial \calU_{j}.
\end{equation}
 The quasi-static approximation ensures that the total volume of condensates is conserved. This follows from integrating equation (\ref{Mullins}a) with respect to $\x\in \Omega\backslash \calU$ and using the divergence theorem:
\begin{equation}
\sum_{j=1}^N \int_{\partial \calU_j}\nabla u(\x)\cdot \n(\x) d\x=0.
\end{equation}
That is, the sum of the fluxes into the $N$ condensates is zero so that there is no net change in the total condensate volume.
The particle flux at the surface $\partial \calU_j$ determines the rate of growth or shrinkage of the droplet under the adiabatic approximation. Introduce a slow time-scale $\tau=t/\epsilon^2 $ and let $r_j(\tau)$ be the slowly varying radius of the $j$-th droplet. The number $\calN_j(\tau)$ of particles within the $j$-th droplet then evolves according to the equation
\begin{align}
\label{calN}
\frac{d\calN_j(\tau)}{d\tau} &= -\epsilon^2 D\int_{\partial \calU_j(\tau)}{\bm \nabla} u_{\tau}(\x) \cdot \n_jd\x ,
\end{align}
where $u_{\tau}(\x)$ is the steady-state concentration for the set of radii $\{r_1(\tau),\ldots,r_N(\tau)\}$, and $\partial \calU_j(\tau)$ is the corresponding droplet boundary.
In the case of a uniform spherical droplet in the high density phase $\phi_b$, the number of particles is $\calN_j(\tau)=|\calU_j(t)|\phi_b$. Hence, equation (\ref{calN}) yields a closed equation for the dynamics of the droplet radii. 

\begin{figure}[t!]
  \centering
  \includegraphics[width=8cm]{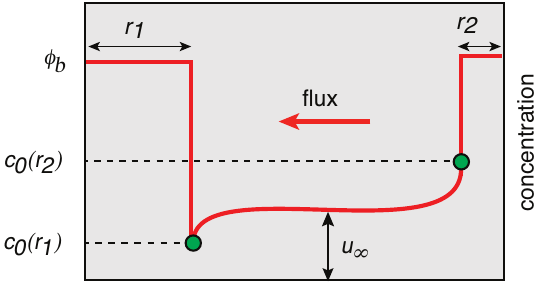}
  \caption{Schematic diagram illustrating the concentration profile as a function of $x$ along the axis joining the centres of two well separated droplets with different radii $r_1>r_2$. The solute concentration $c_0(r_1)$ around the larger droplet is lower than the concentration $c_0(r_2)$ around the smaller droplet due to surface tension. This results in a net diffusive flux from the small droplet to the large droplet. Here $u_{\infty}$ denotes the far-field concentration.}
  \label{fig4}
\end{figure}

In classical Ostwald ripening, the difference in surface concentrations $c_0(r_i)$ for droplets of different sizes results in a net diffusive flux from small to large droplets. This is illustrated in Fig. \ref{fig4} for two droplets $\calU_1$ and $\calU_2$ with $r_2 < r_1$. These fluxes modify the far-field solute concentration $u_{\infty}$, which ultimately leads to complete phase separation\footnote{Standard treatments of Ostwald ripening use a mean-field approximation, in which each droplet ``sees'' an effective far-field concentration $u_{\infty}$ that is determined self-consistently. The theory of diffusion in singularly perturbed domains can be used to calculate corrections to mean field theory \cite{Kavanagh14,Bressloff20a,Bressloff20b,Bressloff24b}.}. However, multiple coexisting biological condensates are found  in the cytoplasm and cell nucleus (see the reviews \cite{Hyman14,Brangwynne15,Berry18,Mittag22} and references therein). These membrane-less organelles are viscous, liquid-like structures containing enhanced concentrations of various proteins and RNA, many of which can be continually exchanged with the surrounding medium. 
 It has been hypothesised that the coexistence of multiple droplets over significant timescales involves the
active suppression of Ostwald ripening. Various hypotheses have been given to account for this suppression: actively driven chemical reactions that maintain the out-of-equilibrium switching of proteins between soluble and phase separating forms \cite{Zwicker15,Wurtz18,Lee19,Weber19} (see also section 6.2);  the mechanical suppression of coarsening mediated by intracellular visco-elastic networks such as the cytoskeleton \cite{Feric13,Style18,Ros20}; the slow conversion of a molecular constituent between the dilute and dense phases due to various regulatory interfacial proteins \cite{Seydoux18,Folk21,Lee21}. The last mechanism can be implemented by considering model I with $\kappa_j<\infty$ \cite{Bressloff24a}.

\medskip

\noindent \underline{Protein clustering in membranes \cite{Bressloff24b}.} A second example of a self-organising process based on model I is protein clustering in cell membranes, where $\calU_j$ is an anchored protein cluster and $u(\x,t)$ is the concentration of free proteins. Each boundary $\calU_j$ is taken to be totally absorbing, and proteins are exchanged between the cytosol and membrane via the active mechanisms of endocytosis and exocytosis \cite{Foret12}. That is, $\gamma_0$ and $\calI_0$ specify the rates of endocytosis and exocytosis, respectively.
The aggregation or clustering of proteins and other macromolecules plays an important role in the formation of various large-scale structures in cell membranes, including the assemblies of cell-cell adhesion proteins in epithelia \cite{Quang13}, lipid raft formation \cite{Turner05,Turner18}, and postsynaptic domains (PSDs) at excitatory and inhibitory synapses in neurons. 
PSDs are rich in scaffolding proteins that can transiently trap transmembrane neurotransmitter receptors, thus localising them at specific sites that are apposed to active zones in presynaptic domains where neurotransmitters are released. PSDs thus play a crucial role in determining the effective strength of synaptic connections between cells \cite{Mac11,Bied17,Choquet18}, see also the next application. One model of PSD formation is based on Smoluchowski coagulation equations, in which the system is kept out of equilibrium by the recycling of particles between the cell membrane and interior \cite{Ranft17,Hakim20}. This results in the stabilisation of a stationary distribution consisting of multiple clusters. (Other models treat PSD formation as a reaction-diffusion process undergoing non-equilibrium pattern formation \cite{Hasel11,Hasel15,Hakim20,Specht21}.) Recently, a modified PSD aggregation model has been introduced, which involves fixed anchoring sites that trap clusters at specific positions within the membrane, consistent with the alignment of PSDs and presynaptic active zone \cite{Ranft23}. These anchoring sites can be identified with the compartmental centres $\x_j$. In this particular model, cluster growth and shrinkage is based on the absorption of
free proteins at the boundaries $\calU_j$, combined with endocytosis within the cluster. There exists a unique steady-state solution satisfying the system of equations
\begin{subequations} 
\label{SScluster}
\begin{align}
	D{\bm \nabla}^2 u(\x)-\gamma_0 u(\x)+\calI_0&=0,\  \x\in \Omega \backslash\calU,\\
	D{\bm \nabla} u(\x) \cdot \n(\x)&=0,\  \x\in\partial \Omega ,\quad
	u(\x)=0,\quad \x\in \partial \calU_j.	\end{align}
	\end{subequations} 
with the radii determined self-consistently from 
\begin{equation}
\label{SScon}
D\int_{\partial \calU_j}{\bm \nabla} u(\x) \cdot \n_jd\x=-\gamma_0 u_0 |\calU_j|,
 \end{equation}
 where $u_0$ is the uniform concentration within a cluster.

\medskip

\noindent \underline{Synaptic receptor trafficking \cite{Bressloff08,Bressloff23a}.}
One application of model II is synaptic receptor trafficking within the cell membrane of a neuron. Advances in single particle tracking (SPT) and imaging techniques have established that freely diffusing neurotransmitter receptors are temporarily confined by binding to scaffolding proteins within the PSDs \cite{Choquet13}. (The dynamics of scaffolding proteins considered in the previous example occurs on slower time scales.) Surface receptors are also exchanged between the cell membrane and interior of the cell via active forms of vesicular transport (exocytosis and endocytosis). The lateral diffusion and trapping of synaptic receptors occurs for almost all types of synapses. This includes the majority of fast excitatory synapses in the central nervous system, which involve the neurotransmitter glutamate binding to $\alpha$-amino-3-hydroxy-5-methyl-4-isoxazole-propionic acid receptors (AMPARs) within dendritic spines \cite{Borgdorff02,Groc04,Thoumine12,Choquet18}, and inhibitory synapses containing the glycine receptor (GlyR) that are found in the postsynaptic membrane of the soma and initial portion of dendrites in spinal cord neurons \cite{Meier01,Dahan03,Triller23}. As the effective strength or weight of a synapse depends on the number of synaptic receptors, it follows that the synaptic weight is determined by a non-equilibrium steady-state in which there is a dynamical balance between different non-zero receptor fluxes (diffusive, endocytotic etc.)  This also implies that activity-dependent changes in the strength of the synapse correspond to shifts in the dynamical balance-point.

Since a typical dendrite is hundreds of microns in length but only a few microns in diameter, most diffusion-trapping models of AMPAR trafficking treat the dendrite as a quasi-one-dimensional (1D) cable along which the spines are represented as discrete point sources or sinks \cite{Earnshaw06,Holcman06,Bressloff07,Earnshaw08,Schumm22}. On the other hand, the reduction to a 1D cable model is not appropriate for inhibitory synapses located in the somatic membrane, for example. In such cases, one has to treat the somatic membrane as a singularly perturbed 2D domain containing one or more synapses that act as transient traps \cite{Bressloff23a}\footnote{Analyzing the full 2D model in the case of a long thin cylindrical surface also establishes the validity of the quasi-1D approximation for dendrites \cite{Bressloff08}.}.
In terms of equations (\ref{gmaster}a,b) and (\ref{gmasterj}a,b), $\calU_j$ represents a synapse, $u(\x,t)$ is the concentration of extrasynaptic receptors and $v_j(\x,t)$ is the concentration of synaptic receptors within the $j$-th synapse. The constants $\gamma_0$ and $\overline{\gamma}_j$ represent the extrasynaptic and synaptic rates of endocytosis, while $\calI_0$ and $\overline{\calI}_j$ are the corresponding rates of exocytosis. The inclusion of a semi-permeable interface
 is motivated by the so-called partitioned fluid-mosaic model of the plasma membrane, in which confinement domains are formed by a fluctuating network of cytoskeletal fence proteins combined with transmembrane picket proteins that act as fence posts \cite{Kusumi05}. Note that the effects of scaffolding proteins within the PSD can be encoded by the synaptic diffusivity $\overline{D}_j$ \cite{Bressloff23a}. This exploits a well-known biophysical mechanism for reducing the diffusion coefficient, in which scaffold proteins act as mobile buffers \cite{Keener09}. Finally, solving the steady-state equations then determines the number of receptors in the $j$-th synapse according to
 \begin{equation}
\label{rk}
\calN_j =\int_{\calU_j}v_j(\x)d\x.
\end{equation}

\begin{figure}[t!]
  \centering
      \includegraphics[width=8cm]{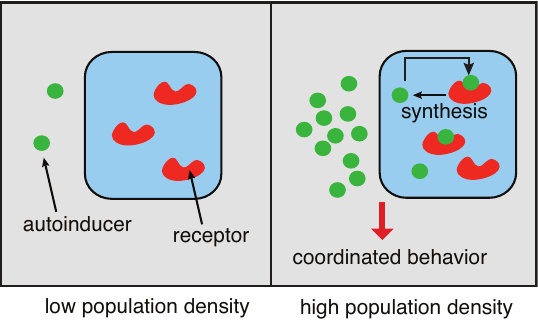}
  \caption{A schematic illustration of quorum sensing at the single-cell level. }\label{fig5}
\end{figure}

\noindent \underline{Bacterial quorum sensing \cite{Muller13,Gou16,Iyan21a,Iyan21b}.} One of the major applications of model III is bacterial quorum sensing (QS), which is a form of collective chemical sensing and response that depends on population density\footnote{A modified version of the PDE-ODE system given by equations (\ref{gmaster}) and (\ref{mass}) has been used to model the spatial spread of an airborne disease between localized patches of human populations \cite{David20}. The dynamics within a patch $\calU_j$ is given by an SIR (susceptible-infected-recovered) model, with the rate of infection determined by the concentration of pathogen on the boundary $\partial \calU_j$. The field $u(\x,t)$ represents the concentration of diffusing pathogen outside the local populations. For simplicity, the movement of people between patches is ignored.}. Examples include  bioluminescence, biofilm formation, virulence, and antibiotic resistance
\cite{Waters05,Hense15,Perez16,Papenfort16,Muk19}. Bacterial QS involves the exchange of certain signalling molecules called autoinducers between the intracellular and extracellular environments. At low cell densities the concentration of autoinducers within individual cells is too low to trigger a cellular response. However, as the population grows, the concentration of autoinducers passes a threshold, resulting in the activation of various genes including those responsible for synthesising the autoinducers. The resulting positive feedback loop means that all of the cells 
initiate transcription at approximately the same time, resulting in some form of coordinated behaviour. The basic process at the single-cell level is shown in Fig. 
\ref{fig5}.

Most models of bacterial QS assume that extracellular diffusion is sufficiently fast so that one can treat the system as a population of identical homogeneous cellular compartments diffusively coupled to a well-mixed extracellular compartment $\Omega\backslash \calU$. The corresponding dynamics is given by a system of coupled nonlinear ODEs of the form
\begin{subequations}
\begin{align}
\label{mass0}
\frac{dU}{dt}&= \frac{\alpha \widehat{\kappa}}{N}\sum_{j=1}^N  \bigg [w_{j,0}(t)-U(t)\bigg ]-\gamma_0 U(t),\quad \alpha =\frac{|\calU|}{|\Omega|-|\calU|},\\
\frac{dw_{j,a}}{dt}&=f_a(\w_j)- \widehat{\kappa}\delta_{a,1} [ w_{j,0}(t)- U(t)] ,
\end{align}
\end{subequations}
where $U(t)$ is the uniform extracellular autoinducer concentration, $w_{j,0}(t)$ is the autoinducer concentration within the $j$-th cell, and $\alpha$ is the volume fraction of the cytosol relative to the extracellular domain. (Note that $\alpha$ is related to the cell density.)  However, in order to take into account the effects of bulk diffusion, it is necessary to consider the coupled PDE-ODE model given by equations (\ref{gmaster}) and (\ref{mass}) with $\gamma_0$ a cytosolic degradation rate, $\calI_0=0$ and $c_j(\x,t)=w_{j,0}(t)$.

 \setcounter{equation}{0}
\section{Matched asymptotics: steady-state analysis}

Consider the steady-state diffusion problem for a constant boundary field $c_{j,0}$ (model I):
\begin{subequations}
\label{masterI}
\begin{align}
 &D{\bm \nabla}^2u(\x)- \gamma_0 {u}(\x) +\calI_0=0,  \quad \x\in \Omega\backslash \calU,\\
&D{\bm \nabla} u(\x) \cdot \n=0,\quad \x \in \partial \Omega,\\
&	D{\bm \nabla} u(\x)\cdot \n_j = \frac{\kappa_j }{\epsilon }[  u(\x) -c_{j,0}] ,\quad \x \in \partial \calU_j.
 \end{align}
\end{subequations}
The corresponding steady-state equations for model II are 
\begin{subequations}
\label{masterII}
\begin{align}
 &D{\bm \nabla}^2u(\x)- \gamma_0 {u}(\x)+\calI_0 =0,  \quad \x\in \Omega\backslash \calU,\\
&D{\bm \nabla} u(\x) \cdot \n=0,\quad \x \in \partial \Omega,\\
&	\overline{D}_j {\bm \nabla}^2 v_j(\x,t)  -\overline{\gamma}_j v_j(\x,t)+\overline{\calI}_j=0, \ \x\in \calU_j, \\
&D{\bm \nabla} u(\x^+)\cdot \n_j(\x^+) = \overline{D}_j{\bm \nabla} v_j(\x^-)\cdot \n_j(\x^-) =\kappa_j[u(\x^+) -v_j(\x^-)] ,\quad \x  \in \partial \calU_j.
 \end{align}
\end{subequations}
Finally, the steady-state version of model III is obtained by setting $c_{j,0}=w^*_{j,0}$ in equation (\ref{masterI}),with
\begin{align}
\label{masterIII}
f_a(\w^*_j)=  \frac{ \kappa_j \delta_{a,1} }{|\calU_j|} \int_{\partial \calU_j}( w^*_{j,0}- u(\x))d\x.
\end{align}

\begin{figure}[t!]
  \centering
  \includegraphics[width=12cm]{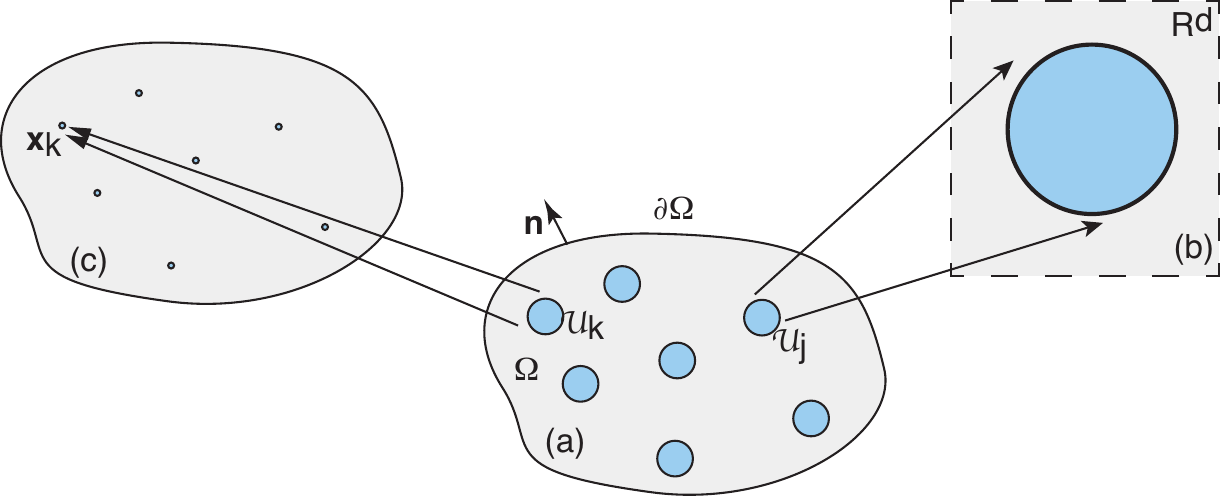}
  \caption{Schematic illustration showing the mapping of a singularly perturbed diffusion problem to corresponding inner and outer solutions. (a) Original unscaled domain. (b) Construction of the inner solution in terms of stretched coordinates $\y=\epsilon^{-1}(\x-{\x}_j)$, where ${\x}_j$ is the centre of the $j$-th compartment. The rescaled radius is $\ell_j$ and the region outside the compartment is taken to be $\R^d$, $d=2,3$, rather than the bounded domain $\Omega$. (c) Construction of the outer solution. Each compartment is shrunk to a single point. The outer solution can be expressed in terms of the corresponding modified Neumann Green's function and then matched with the inner solution around each compartment.}
  \label{fig6}
\end{figure}

Given the basic assumption that the compartments are small and well-separated, we can use a combination of matched asymptotics and Green's function methods along the lines of Refs. \cite{Ward93,Ward93a,Straube07,Gou16,Bressloff23a} in 2D and Refs. \cite{Straube09,Cheviakov11,Bressloff23b} in 3D. More specifically, we construct an inner or local solution
valid in an $O(\epsilon)$ neighbourhood of each compartment, and then match to an outer or global solution that is valid away from each neighbourhood. The general construction is illustrated in Fig. \ref{fig6}. We now give a detailed description of the steps of the analysis, first in 2D and then 3D. In both cases, we develop the asymptotic analysis for model I and then indicate how to extend the analysis models II and III. We also perform the rescalings
\begin{equation}
\kappa_j\rightarrow  \frac{\kappa_j}{\epsilon},\quad  \overline{\calI}_j\rightarrow \frac{\overline{\calI}_j}{\epsilon^2}, \quad \bar{\gamma}_0\rightarrow \frac{\bar{\gamma}_j}{\epsilon^2},\quad j=1,\ldots,N,
\end{equation}
so that the total flux across the boundary $\partial \calU_j$ and the rates of insertion and removal within $\calU_j$ are $O(1)$\footnote{One example where the reactivity $\kappa_j$ is not rescaled is in a conversion-limited model of Ostwald ripening during the liquid-liquid phase separation of biological condensates \cite{Lee21,Bressloff24b} Taking the flux across each droplet boundary to be $O(\epsilon)$ means that the coarsening of droplets is significantly slowed down, resulting in the coexistence of multiple droplets over relevant biological timescales. The corresponding matched asymptotic analysis has to be modified accordingly.}.
Finally, without loss of generality, se take $\calI_0=0$. (If $\calI_0>0$ then we simply perform the shifts $u(\x)\rightarrow u(\x)-\calI_0/\gamma_0$, $c_{j,0}\rightarrow c_{j,0}-\calI_0/\gamma_0$ etc.)

\subsection{Asymptotic analysis in 2D}
In the 2D case ($d=2$), the inner solution of equations (\ref{masterI}) near the $j$-th compartment is constructed by introducing the stretched local variable ${\mathbf y} =
\varepsilon^{-1}(\x-\x_j)$ and setting
\begin{equation}
U_j(\y)=u(\x_j+\varepsilon \y).
\end{equation}
The inner solution $U_j$ then satisfies (on dropping $O(\epsilon)$ terms)
\begin{subequations}
\label{inner}
\begin{align}
&D {\bm \nabla}^2_{\y}U_j(\y)= 0,\quad |\y| > \ell_j,\\
	&D{\bm \nabla}_{\y} U_j(\y)\cdot \n_j(\y)= \kappa_j [U_j(\y)-c_{j,0}],\quad   |\y |=\ell_j  .
\end{align}
\end{subequations}
Using polar coordinates with $|\y|=\rho$ and setting $U_j=U_j(\rho)$, we have
\begin{equation}
\frac{d^2  U_j}{d\rho^2}+\frac{1}{\rho} \frac{d  U_j}{d\rho}=0, \quad \rho>\ell_j,\quad D\left .\frac{d  U_j}{d\rho}\right |_{\rho=\ell_j}=\kappa_j[U_j(\ell_j)-c_{j,0}],
\end{equation}
which has the general solution
\begin{align}
\label{UV}
    &U_j(\rho) =  \Phi_j +\calA_j \log[\rho/\ell_j ],\quad \ell_j \leq \rho < \infty ,
\end{align}
with
\begin{align}
\label{UV2}
 \Phi_j=c_{j,0} +\frac{\calA_jD}{\kappa_j\ell_j}.
 \end{align}
In order to determine the $N$ coefficients $\calA_j$ we have to match the far-field behaviour of the inner solution $U$ with the outer solution. The outer solution is constructed by shrinking each compartment to a single point and imposing a corresponding singularity condition that is obtained by matching with the inner solution. The outer equation is given by
\begin{align}
\label{outer}
	 D{\bm \nabla}^2 u(\x) - \gamma_0 u(\x)= 0
		\end{align}
for $\x\in \Omega'\equiv \Omega\backslash \{\x_1,\ldots,\x_N\},$ together with the boundary condition $D\nabla u \cdot \n=0, \ \x\in \partial \Omega$. The corresponding singularity conditions are
\begin{align}
\label{goo}
   u(\x) \sim   \Phi_j+  \calA_j \left[\log(|{\bf x} - {\bf x}_j|/\ell_j)-\log \epsilon \right]
\end{align}
for ${\bf x} \to {\bf x}_j$. 
A common feature of strongly localised perturbations in 2D domains \cite{Ward93} is the appearance of the small parameter
\begin{equation}
\nu=-\frac{1}{\log \epsilon }.
\end{equation}
In order to eliminate the $1/\nu$ term in the singularity condition (\ref{goo}), we rescale the unknown coefficients $\calA_j$ by setting $\calA_j=\nu A_j $.
 It is well-known that $\nu \rightarrow 0$ more slowly than $\epsilon\rightarrow 0$. Hence, if one is interested in obtaining $O(1)$ accuracy with respect to an $\epsilon$ expansion, then it is necessary to sum over the logarithmic terms non-perturbatively \cite{Ward93,Ward93a}.

The first step is to introduce the Neumann Green's function of the modified Helmholtz equation according to\footnote{There are very few geometries for which the 2D or 3D Neumann Green's function of the modified Helmholtz equation is known explicitly, see appendix A.}
\begin{subequations}
\label{GMH}
\begin{align}
D\nabla^2 G(\x,\x') - \gamma_0 G(\x,\y)&=-\delta(\x-\x') , \quad \x,\x'\in \Omega,\\
	\nabla G(\x,\x')\cdot \n(\x)&=0 ,\quad \x \in \partial \Omega,\quad \int_{\Omega} G(\x,\x')d\x=\frac{1}{\gamma_0} . 
	\end{align}
	\end{subequations}
Note that the 2D Green's function can be decomposed as 
\begin{align}
    G(\x,\x') = -\frac{\log{|\x - \x'|}}{2\pi D} + R(\x, \x'),
\end{align}
where $R$ is the non-singular part of the Green's function. It follows that the solution of equation (\ref{outer}) can be written in the form 
\begin{align}
    u(\x) &\sim   - 2\pi \nu D\sum_{k=1}^N A_k G(\x,\x_k).
    \label{outu}
\end{align}
We have $N$ unknown coefficients $A_j $, which are determined by matching the inner and outer solutions:
\begin{align}
   A_j \left[1 + \nu \Psi_j-\nu \log \ell_j + 2\pi  D \nu R(\x_j,\x_j )\right]+ 2\pi  D \nu \sum_{k \neq j}A_k G(\x_j, \x_k) =-c_{j,0}.
  \label{match}
\end{align}
We have set $\Phi_j=c_{j,0}+\nu A_j(\nu)\Psi_j$ for convenience, with $\Psi_j=D/(\kappa_j\ell_j)$.
Let us rewrite equation (\ref{match}) as a matrix equation:
\begin{align}
\label{matrix}
    \left[{\bf I} + \nu\left({\bm \Psi} + 2\pi  D{\bm \calG} \right)\right]{\bf a} = -{\bf c}_0,
\end{align}
where ${\bf I}$ is the $N \times N$ identity matrix, 
\begin{align}
{\bf a}=(A_1, \ldots, A_N)^{\top},\quad {\bm \Psi}=\mbox{diag}(\Psi_1,\ldots,\Psi_N), 
\end{align}
and ${\bm \calG}$ is an $N \times N$ matrix with entries
\begin{equation}
\calG_{ij}=G(\x_i,\x_j),\ i\neq j,\quad \calG_{jj}=R(\x_j,\x_j)-\frac{\log \ell_j}{2\pi D}.
\end{equation}
Inverting equation (\ref{matrix}) yields in component form
\begin{align}
\label{naj}
    A_j 
    &= -\sum_{i=1}^N\left({\bf I}  + \nu {\bf M} \right)^{-1}_{ji}c_{i,0} ,\quad j = 1, \ldots, N,
    \end{align}
 with
  \begin{equation}
  \label{M}
 M_{ji}=\Psi_j\delta_{i,j} +2\pi  D \calG_{ji}.
    \end{equation}
Note that the solution (\ref{naj}) for $A_j$ is a non-perturbative function of the small parameter $\nu$, which was obtained by matching the inner and outer solutions using Green's functions along the lines originally developed in Refs. \cite{Ward93,Ward93a}. This effectively sums over the logarithmic terms, which is equivalent to calculating the asymptotic solution for all terms of $O(\nu^k)$ for any $k$. It is $O(1)$ with respect to a corresponding $\epsilon$ expansion. If one were to Taylor expand equation (\ref{naj}) in powers of $\nu$, then the leading order form of the outer solution would be
 \begin{align}
    u(\x) &\sim    2\pi \nu D\sum_{k=1}^N c_{k,0}G(\x,\x_k).
\end{align}
 Finally, note that the asymptotic solution (\ref{naj}) still holds in the limit $\kappa_j\rightarrow \infty$, except that $\Psi_j\rightarrow 0$. This corresponds to having the Dirichlet boundary condition $u(\x)=c_{j,0}$ for all $\x\in \partial \calU_j$.
 
 \paragraph{\underline{Zero degradation ($\gamma_0=0$).}} Special care has to be taken when $\gamma_0=0$, since the solution of Laplace's equation in a bounded domain $\Omega$ with a reflecting boundary $\partial \Omega$ is only defined up to a constant. The corresponding generalised Neumann Green's function satisfies the equation 
 \begin{subequations}
\label{exG1}
\begin{align}
D{\bm \nabla}^2 G_0&=\frac{1}{|\Omega|}-\delta(\x-\x'),\quad \x,\x' \in \Omega,\\
\n \cdot {\bm \nabla} G_0& =0 \mbox{ on } \partial \Omega,\quad  \int_{\Omega}G_0 d\x=0
\end{align}
\end{subequations}
for fixed $\y$. The additional constant term $1/|\Omega|$ on the right-hand side of equation (\ref{exG1}a) is needed in order to ensure that both sides yield zero when integrating with respect to $\x \in \Omega$. The condition $\int_{\Omega}G_0 d\x=0$ determines $G_0$ uniquely. Again we can decompose $G_0 $ as
\begin{equation}
G_0(\x,\x')=-\frac{ \ln |\x-\x'|}{2\pi D}+R_0(\x,\x'),
\end{equation}
where $R_0$ is the regular part of the Green's function $G_0$. (Some examples of $G_0$ in simple geometries are collected in appendix A.) The analysis of the inner solution is identical to the case $\gamma_0>0$. However, the outer solution is now only defined up to a constant $u_{\infty}$ such that
\begin{equation}
\label{exouter}
u(\x)\sim  u_{\infty}-2\pi  D \nu \sum_{k=1}^NA_k G_0(\x,\x_k)
\end{equation}
for $\x \notin \{\x_j,\, j=1,\ldots,N\}$
and some constant $u_{\infty}$. Since $\int_{\Omega} G_0d\x=0$, 
it follows that
\begin{equation}
u_{\infty}=|\Omega|^{-1}\int_{\Omega}u(\x)d\x.
\end{equation}
In addition, if $\x \notin \{\x_j,\, j=1,\ldots,N\} $ then
\begin{align*}
{\bm \nabla}^2 u(\x) &\sim -2\pi D \sum_{k=1}^NA_k {\bm \nabla}^2 G_0(\x,\x_k) 
= -\frac{2\pi}{|\Omega|}  \sum_{=1}^NA_k.
\end{align*}
Hence, the outer solution satisfies the steady-state diffusion equation if and only if
\begin{equation}
\label{exsumA}
\sum_{k=1}^NA_k=0.
\end{equation}
In order to determine the $N+1$ unknown coefficients, $A_j$, $j=1,\ldots,N$ and $u_{\infty}$, we require $N+1$ linearly independent conditions. One of these is given by equation (\ref{exsumA}) whereas the others are obtained from matching the inner and outer solutions. Proceeding along similar lines to the derivation of equation (\ref{naj}) we find that
\begin{align}
\label{naj0}
    A_j
    &= \sum_{i=1}^N\left({\bf I}  + \nu {\bf M}_0 \right)^{-1}_{ji}[u_{\infty}-c_{i,0}] ,\quad j = 1, \ldots, N,  
    \end{align}
 with
\begin{align}
{\bf M}_0={\bm \Psi}+2\pi D{\bm \calG}_0,\quad \calG_{0,jj}&=R_0 (\x_j,\x_j)-\frac{\ln \ell_j}{2\pi D} \mbox{ and }  \calG_{0,ji}=G_0(\x_j,\x_i),\quad j\neq i.
\label{exM}
\end{align}
Imposing the constraint (\ref{exsumA}) on equation (\ref{naj0}) then determines the unknown constant $u_{\infty}$: 
\begin{align}
\label{exuinf}
u_{\infty}&=\left [
{\displaystyle \sum_{i,j=1}^N \left({\bf I}  + \nu {\bf M}_0 \right)^{-1}_{ji}}\right ]^{-1} {\displaystyle \sum_{i,j=1}^N  \left({\bf I}  + \nu {\bf M}_0 \right)^{-1}_{ji}c_{j,0}} .\end{align}
Finally, note that expanding equations (\ref{naj0}) and (\ref{exuinf}) to $O(\nu)$ yields
\begin{equation}
A_j\sim u_{\infty}-c_{j,0},\quad  u_{\infty}\sim \overline{c}_0\equiv \frac{1}{N} \sum_{k=1}^N c_{k,0},
\end{equation}
and the outer solution reduces to 
\begin{equation}
u(\x)\sim \overline{c}_0-2\pi D\nu \left [\sum_{k=1}^N(\overline{c}_0-c_{k,0})G_0(\x,\x_k)-\frac{1}{N}\sum_{j,k=1}^N (\overline{c}_0-c_{k,0})[G_{0,jk}+\Psi_k\delta_{j,k}]\right ].
\end{equation}

One important application of the above analysis is determining the coarsening dynamics of 2D Ostwald ripening \cite{Kavanagh14,Bressloff20a}. First, substituting the inner solution (\ref{UV}) into equation (\ref{calN}) under the adiabatic approximation and setting  $|\calU_j| =\pi r_j^2$ with $r_j=\epsilon \ell_j$ gives to leading order
\begin{equation}
\frac{d\ell_j}{d\tau} =\frac{D \nu A_j}{\ell_j\phi_b } \sim \frac{D \nu }{\ell_j\phi_b } (u_{\infty}-c_{j,0}),\quad j=1,\ldots,N .
\end{equation}
Substituting for $c_{j,0}$ using equation (\ref{GT0}) then implies that
\begin{equation}
u_{\infty}-c_{j,0} \sim \frac{\phi_a }{N} \sum_{k=1}^N\left (1+\frac{\ell_c}{r_k}\right )-\phi_a\left (1+\frac{\ell_c}{r_j}\right )= \phi_a\left [\frac{1}{N}\sum_{k=1}^N\frac{1}{r_k}-\frac{1}{r_j}\right ].
\end{equation}
Hence,
\begin{equation}
\frac{d\ell_j}{d\tau} =\frac{D \nu \phi_a}{\ell_j\phi_b }\left (\frac{1}{\ell_{\rm harm}}-\frac{1}{\ell_j}\right ),\quad j=1,\ldots,N ,
\end{equation}
where $\ell_{\rm harm}$ is the harmonic mean of the droplet radii,
\begin{equation}
\frac{1}{\ell_{\rm harm}}=\frac{1}{N}\sum_{k=1}^N\frac{1}{r_k}.
\end{equation}
Including higher-order terms in the asymptotic expansion then yields corrections to this mean field result.

\paragraph{\it \underline{Semi-permeable interface (model II)}.}

The only modification of the inner solution (\ref{UV}) and the outer solution (\ref{outu}) in the case of a semipermeable interface concerns the calculation of the constants $c_{j,0}$ and $\Psi_j$ appearing in the matching condition (\ref{match}). That is, we have to extend the inner solution so that it includes the effects of diffusion within each compartment as determined by equations (\ref{masterII}). In terms of stretched coordinates, we set
\begin{equation}
U_j(\y)=u(\x_j+\epsilon \y),\quad V_j(\y)= v_j(\x_j+\epsilon \y).
\end{equation}
The inner solutions $U_j,V_j$ then satisfy (on dropping $O(\epsilon)$ terms)
\begin{subequations}
\label{iperm}
\begin{align}
&D {\bm \nabla}^2_{\y}U_j(\y)= 0,\ |\y| > \ell_j,\\ &\overline{D}_j{\bm \nabla}^2_{\y}V_j (\y)- \overline{\gamma}_j V_j(\y)+\overline{\calI}_j= 0 , \   |\y| < \ell_j ,\\
	&D{\bm \nabla}_{\y} U_j(\y^+)\cdot \n_j(\y^+) = \overline{D}_j{\bm \nabla}_{\y} V_j(\y^-)\cdot \n_j(\y^-) =\kappa_j [U_j(\y^+)- V_j(\y^-)],\quad   |\y  |=\ell_j  .
\end{align}
\end{subequations}
Using polar coordinates with $|\y|=\rho$, the solution can be written as \cite{Bressloff23a}
\begin{subequations}
\label{UVj}
\begin{align}
    &U_j=  \Phi_j +\nu A_j(\nu) \log\rho/\ell_j ,\quad \ell_j \leq \rho < \infty, \\
    &V_j=\frac{\overline{\calI}_j}{\overline{\gamma}_j}+\overline{\Phi}_j \frac{I_0\left(\overline{\beta}_j \rho\right)}{I_0(\overline{\beta}_j  \ell_j)},\quad 0 \leq \rho \leq \ell_j ,
\end{align}
\end{subequations}
where $\overline{\beta}_j=\sqrt{\bar{\gamma}_j/\bar{D}_j}$ and $I_0$ is a modified Bessel function of the first kind.
The coefficients $\Phi_j,\overline{\Phi}_j$ can be expressed in terms of $ \nu A_j$ by imposing equations (\ref{iperm}c). After some algebra, we find that 
\begin{align}
\Phi_j =\frac{\overline{\calI}_j}{\overline{\gamma}_j}+\nu A_j \left [\frac{D}{\kappa_j\ell_j}+ \frac{ D }{\overline{D}_j{\mathcal F}(\beta_j \ell_j)}\right ] ,\quad \overline{\Phi}_j &=\nu A_j\frac{D}{\overline{D}{\mathcal F}(\overline{\beta}_j  \ell_j)},
\label{Psip}
\end{align}
with 
\begin{equation}
{\mathcal F}(x):=\frac{xI_1\left(x\right)}{I_0(x)}.
\end{equation}
It follows that 
\begin{equation}
c_{j.0}= \frac{\overline{\calI}_j}{\overline{\gamma}_j}, \quad \Psi_j=\frac{D}{\kappa_j\ell_j}+ \frac{ D }{\overline{D}_j{\mathcal F}(\overline{\beta}_j \ell_j)}
\end{equation}
in the matching condition (\ref{match}).
Note that in the limit $\overline{D}_j\rightarrow \infty$ for all $j=1,\ldots,N$ we recover model I.

The above model has been used to study the effects of bulk diffusion on synaptic receptor trafficking in somatic inhibitory synapses \cite{Bressloff23a}.
Substituting equation (\ref{UVj}b) into equation (\ref{rk}) then gives the following expression for the steady-state number of receptors in the $j$-th synapse: 
 \begin{align}
\calN_j&= 2\pi \epsilon^2 \int_{\calU_j}V_j(\rho)\rho d\rho\sim  \pi \epsilon^2\ell_j^2  \frac{\overline{\calI}_j}{\overline{\gamma}_j}+ \frac{2\pi \epsilon^2 \nu DA_j }{D_j}\int_0^{\ell_j}\frac{I_0\left(\overline{\beta}_j\rho \right)}{\overline{\beta}_j \ell_jI_1(\overline{\beta}_j\ell_j)}\rho d\rho\nonumber \\
&=\pi \epsilon^2\ell_j^2  \frac{\overline{\calI}_j}{\overline{\gamma}_j}+\frac{2\pi \nu DA_j }{\overline{\gamma}_j/\epsilon^2}.
\label{rstar}
\end{align}
We have used the Bessel function identity $\frac{d}{dx}[xI_1(x)]=xI_0(x)$, which implies that
\begin{align*}
\int_0^{\ell_j} I_0(\beta_j\rho)\rho d\rho&=\beta_j^{-2} \int_0^{\beta_j \ell_j}x I_0(x)dx
=\frac{1}{\overline{\beta}_j^{2}} \int_0^{\overline{\beta}_j\ell_j}\frac{d}{dx} [x I_1(x)]dx=\frac{\ell_j}{\overline{\beta}_j}I_1(\overline{\beta}_j \ell_j).
\end{align*}
It follows that the synaptic strength as determined by $\calN_j$ depends on local synaptic parameters such as the rates of endo/exocytosis $\overline{\gamma}_j/\epsilon^2$ and $\overline{\calI}_j/\epsilon^2$, the synaptic diffusivity $\overline{D}_j$, and the size of the synapse $\epsilon \ell_j$. However, there is also a heterosynaptic component to the synaptic strength due to the fact that the coefficient $A_j$ depends on all of the other synapses through the matching condition (\ref{match}) \cite{Bressloff23a}. 

\paragraph{\underline{Biochemical reaction network (model III).}}

In compartmental model III the constant $c_{j,0}=w_{j,0}^*$ is determined from equation (\ref{masterIII}) under the mappings $\kappa_j\rightarrow \kappa_j/\epsilon$ by setting $u(\x)=U_j(\ell_j)$ on $\partial \calU_j$:
\begin{align}
\label{stmass}
\widehat{f}_a(\w_j^*)=   \kappa_j \delta_{a,1}  \int_{\partial \calU_j}( w^*_{j,0}- U_j(\ell_j))d\y= \kappa_j \delta_{a,1} 2\pi \ell_j(w_{j,0}^*- U_j(\ell_j)).
\end{align}
We have also set $f_a(\w)=\widehat{f}_a(\w)/|\calU_j|$.
Using equations (\ref{UV}) and (\ref{UV2}) with $c_{j,0}=w_{j,0}^*$, we obtain the nonlinear algebraic system \cite{Gou16}
\begin{equation}
\label{Fj}
\widehat{f}_a(\w_j^*)+ 2\pi D\nu A_j \delta_{a,1}=0.
\end{equation}
(Recall from equation (\ref{naj}) that the coefficients $A_j$ depend on all the coefficients $c_{k,0}$, $k=1,\ldots,N$.) Since the underlying coupled PDE-ODE system is nonlinear, it follows that there may exist more than one solution of equation (\ref{Fj}). However, in order to determine the stability of each steady-state solution, one has to linearise the full time-dependent equations and solve a non-trivial eigenvalue problem \cite{Gou16}. Linear stability analysis also yields necessary conditions for the occurrence of a Hopf bifurcation, signalling the emergence of collective oscillations. Alternatively, asymptotic methods can be used to reduce the PDE-ODE system to an effective nonlinear ODE system in the fast diffusion regime \cite{Gou16,Iyan21a,Iyan21b}, providing a more tractable framework for investigating diffusion-mediated synchronisation, see also section 5.

\subsection{Asymptotic analysis in 3D}

The 3D analysis of equations (\ref{masterI}) is based on an asymptotic expansion with respect to $\epsilon$ rather than $\nu$ \cite{Straube09,Cheviakov11,Bressloff23b}. The outer solution is expanded as
\begin{equation}
u(\x)\sim  \epsilon u_1(\x)+\epsilon^2 u_2(\x)+\ldots
\end{equation}
with
\begin{subequations}
\label{3Douter}
\begin{align}
D{\bm \nabla}^2 u_{n}(\x) -\gamma_0 u_{n}(\x)=0,\quad \x\in \Omega\backslash \{\x_1,\ldots,\x_N\}; \quad {\bm \nabla}u_{n}(\x)\cdot \n(\x)=0,\quad \x \in \partial \Omega
\end{align}
\end{subequations}
for $n=1,2.\ldots$. Equations (\ref{3Douter}) are supplemented
by singularity conditions as $\x\rightarrow \x_j$, $j=1,\ldots,N$. The latter are determined by matching to the inner solution. In the inner region around the $j$-th compartment, we again introduce the stretched coordinates ${\bf y}=\epsilon^{-1}(\x-\x_j)$ and set $U_j({\bf y}) =u(\x_j+\epsilon \y)$. Expanding the inner solution as 
\begin{equation}
U_j(\y) = U_{j,0}(\y) +\epsilon U_{j,1}(\y)+\ldots,
\end{equation}
we have
\begin{subequations}
\label{3Dinner0}
\begin{align}
&D {\bm \nabla}^2_{\y}U_{j,n}(\y)= \gamma_0 U_{j.n-2}(\y),\quad |\y| > \ell_j,\\
	&D{\bm \nabla}_{\y} U_{j,n}(\y)\cdot \n_j(\y) - \kappa_j U_{j,n}(\y)=-\kappa_j c_{j,0}\delta_{n,0},\quad   |\y |=\ell_j  .
\end{align}
\end{subequations}
(Note that $U_{j,n}\equiv 0$ for $n=-2,-1$.) 
Finally, the matching condition is that the near-field behaviour of the outer solution as $\x\rightarrow \x_j$ should agree with the far-field behaviour of the inner solution as $|\y|\rightarrow \infty$, which is expressed as 
\begin{equation}
\label{3Dmatch}
  \epsilon u_1(\x)+\epsilon^2 u_2(\x) \sim U_{j,0}(\x)+  \epsilon U_{j,1}(\x)+\ldots
\end{equation}
In particular, note that the far-field behaviour of $U_{j,n}$ determines the near-field behaviour of $u_{n+1}$ so we alternate between the inner and outer solutions during matching.

In order to solve the inner BVP for $U_{j,0}$, we introduce spherical polar coordinates and set $U_{j,0}=U_{j,0}(\rho)$: 
\begin{equation}
\frac{d^2  U_{j,0}}{d\rho^2}+\frac{2}{\rho} \frac{d  U_{j,0}}{d\rho}=0, \quad \rho>\ell_j,\quad D\left .\frac{d  U_{j,0}}{d\rho}\right |_{\rho=\ell_j}=\kappa_j[U_{j,0}(\ell_j)-c_{j,0}].
\end{equation}
Assuming that the solution decays at infinity, we have
\begin{equation}
\label{Lam}
 U_{j,0}=\frac{\Lambda_jc_{j,0} }{\rho},\quad \Lambda_j:=\frac{\kappa_j \ell_j^2 }{\kappa_j\ell_j+D}.
\end{equation}
(Note that $\Lambda_j\rightarrow \ell_j$ in the limit $\kappa_j\rightarrow \infty$.)
The matching condition (\ref{3Dmatch}) then implies that $u_1(\x)$ satisfies equations (\ref{3Douter}) together with the singularity condition
\[u_1(\x)\sim \frac{\Lambda_jc_{j,0}}{|\x-\x_j|} \quad \mbox{as } \x\rightarrow \x_j.\]
It follows that the leading order contribution to the outer solution can be written as
\begin{equation}
u_1(\x)=4\pi D \sum_{k=1}^N\Lambda_kc_k G(\x,\x_k),
\end{equation}
where $G$ is the 3D version of the modified Helmholtz Green's function defined in equation (\ref{GMH}). In particular, the 3D Green's function has the singularity structure
\begin{equation}
G(\x,\x')=\frac{1}{4\pi D|\x-\x'|} +R(\x,\x').
\end{equation}
Next we match the far-field behaviour of $U_{j,1}(\y)$ with the non-singular contribution to the near-field behaviour of $u_1(\x)$ around $\calU_j$. That is, $U_{j,1}(\y)\sim \chi_j $ as $ |\y|\rightarrow \infty$ with
\begin{align}
\label{chi}
\quad \chi_j=4\pi D  \Lambda_j c_{j,0}R(\x_j,\x_j)+ 4\pi D \sum_{k\neq j}^N\Lambda_k c_{k,0}G(\x_j,\x_k).
\end{align}
It follows that the solution of equation (\ref{3Dinner0}) for $n=1$ is 
\begin{equation}
U_{j,1} =\chi_ j\left ( 1-\frac{\Lambda_j}{\rho}\right).  
\end{equation}
Hence, $u_2(\x)$ satisfies equation (\ref{3Douter}) supplemented by the singularity condition
\[u_2(\x)\sim - \frac{\chi_j\Lambda_j }{|\x-\x_j|} \quad \mbox{as } \x\rightarrow \x_j.\]
Following along identical lines to the derivation of $u_1(\x)$, we obtain the result
\begin{equation}
u_2(\x)=-4\pi D\sum_{k=1}^N\Lambda_k \chi_k G(\x,\x_k).
\end{equation}
In conclusion, the outer solution takes the form
\begin{align}
\label{out3D0}
u(\x) 
 \sim {4\pi}  D \epsilon \sum_{k=1}^N\Lambda_k(c_{k,0}-\epsilon \chi_k)G(\x,\x_k) , 
\end{align}
while the inner solution around the $j$-th cluster is
 \begin{align}
 \label{Ujr}
U_j(\rho)
\sim  \frac{\Lambda_j}{\rho} \left (c_{j,0}-\epsilon \chi_j\right ).
\end{align}

\paragraph{\underline{Zero degradation ($\gamma_0=0$).}} As in the case of 2D diffusion, we have to modify the analysis when $\gamma_0=0$. In particular, the outer solution has a constant $O(1)$ term that has to be determined self-consistently:
\begin{equation}
u(\x)\sim u_{\infty}+\epsilon u_1(\x)+\epsilon^2 u_2(\x)+\ldots
\end{equation}
The leading order contribution to the inner solution now takes the form
\begin{equation}
 U_{j,0}(\rho) =u_{\infty} -\frac{\Lambda_j(u_{\infty}-c_{j,0})}{\rho},
 \end{equation}
 with $\Lambda_j$ defined in (\ref{Lam}). The subsequent analysis proceeds along identical lines to the previous case. The outer solution becomes
 \begin{align}
 \label{out3D}
u(\x) 
 \sim u_{\infty} +{4\pi}  D \epsilon \sum_{k=1}^N\Lambda_k(c_{k,0}-u_{\infty} -\epsilon \chi_k)G_0(\x,\x_k) , 
\end{align}
while the inner solution around the $j$-th cluster is now
 \begin{align}
U_j(\rho)
\sim  u_{\infty} +\frac{\Lambda_j}{\rho} \left (c_{j,0}-u_{\infty} -\epsilon \chi_j\right ). 
\end{align}
Here $G_0$ is the 3D version of the modified Green's function of equation (\ref{exG1}), and $\chi_j$ is given by equation (\ref{chi}) after replacing $c_{j,0}$ by $c_{j,0}-u_{\infty}$. The unknown constant $u_{\infty}$ is determined by ensuring that the outer solution (\ref{out3D}) satisfies the steady-state equation ${\bm \nabla}^2u(\x)=0$. Applying the Laplacian ${\bm \nabla}^2$ to both sides of equation (\ref{out3D}) implies that for $\x \notin \{\x_j,\, j=1,\ldots,N\} $, we have
\begin{align*}
0=\nabla^2 u(\x) & \sim 4\pi D \sum_{j=1}^N \Lambda_j \left (c_{j,0}-u_{\infty} -\epsilon \chi_j\right ){\bm \nabla}^2 G_0(\x,\x_i)
\\
&\sim \frac{4\pi}{|\Omega|}  \sum_{j=1}^N \Lambda_j \left (c_{j,0}-u_{\infty} -\epsilon \chi_j\right ).
\end{align*}
It follows that to $O(\epsilon)$
\begin{equation}
\label{3Dcon}
u_{\infty}=\frac{ \sum_{j=1}^N\Lambda_j[c_{j,0}-\epsilon \chi_j]}{ \sum_{i=1}^N\Lambda_i }.
\end{equation}

Equations of the form (\ref{out3D0}) and (\ref{out3D}) have been used to characterise the spatial decay of intracellular protein gradients under phosphorylation- dephosphorylation \cite{Straube09}. 
Analogous to the 2D case, equation (\ref{out3D}) has also been used to determine droplet coarsening dynamics in 3D Ostwald ripening.
Substituting equation (\ref{GT0}) into equation (\ref{3Dcon}) and taking the limit $\Lambda_j\rightarrow \ell_j$, we have
\begin{equation}
\label{3Dcon2}
u_{\infty}\sim \frac{\phi_a \sum_{j=1}^N\ell_j \left (1+{\ell_c}/{\epsilon \ell_j}\right )}{ \sum_{i=1}^N\ell_i }=\phi_a+\frac{\phi_a\ell_c}{\sum_{i=1}^N\ell_i}.
\end{equation}
Moreover, since $|\calU_j|=4\pi r_j^3/3$ with $r_j=\epsilon \ell_j$, we find that equation (\ref{calN}) becomes
\begin{equation}
\frac{d\ell_j}{dt} =\frac{D\phi_a\ell_c}{ \phi_b\ell_j }\left (\frac{1}{\ell_{\rm av}}-\frac{1}{\ell_j}\right ),\quad j=1,\ldots,N,\quad \mbox{where } \ell_{\rm av}=\frac{1}{N}\sum_{k=1}^N\ell_k.
\label{radOR}
\end{equation}
Hence, at time $\tau$ all droplets satisfying $\ell_j(\tau) <\ell_{\rm av}(\tau)$ ($\ell_j(\tau) <\ell_{\rm av}(\tau)$) are growing (shrinking). This means that the critical radius $\ell_{\rm av}(\tau)$ is an increasing function of time so that only a single droplet remains in the limit $\tau \rightarrow \infty$.

\paragraph{\underline{Semi-permeable interface (model II).}} Consider the 3D version of the inner BVP (\ref{iperm}). Introducing the series expansions 
\begin{equation}
U_j(\y) = U_{j,0}(\y) +\epsilon U_{j,1}(\y)+\ldots,\quad V_j(\y) = V_{j,0}(\y) +\epsilon V_{j,1}(\y)+\ldots
\end{equation}
we have for $n=0,1$
\begin{subequations}
\label{3Dinner}
\begin{align}
&D {\bm \nabla}^2_{\y}U_{j,n}(\y)=0,\\ &  \overline{D}_j {\bm \nabla}^2_{\y}V_{j,n}(\y)-\overline{\gamma}_j V_{j,n}(\y)+\overline{I}_j\delta_{n,0}=0,\quad |\y| < \ell_j,\\
& D{\bm \nabla} U_{j,n}(\y) \cdot \n_j  =   \overline{D}_j{\bm \nabla} V_{j,n}(\y) \cdot \n_j
 = \kappa_j  [U_{j,n}(\y) -V_{j,n}(\y) ] ,\quad |\y|=\ell_j.
\end{align}
\end{subequations}
For the sake of illustration, we calculate the leading order contribution ($n=0$), which in spherical polar coordinates is the solution to the BVP
\begin{subequations}
\label{sph}
\begin{align}
 &D\frac{\partial^2U_{j,0}(\rho)}{\partial \rho^2} + \frac{2D}{\rho}\frac{\partial U_{j,0}(\rho)}{\partial \rho}=0, \ \rho>\ell_j,\\
&\overline{D}_j\frac{\partial^2V_{j,n0}(\rho)}{\partial \rho^2} + \frac{2\overline{D}}{\rho}\frac{\partial V_{j,0}(\rho)}{\partial \rho}-\overline{\gamma}_jV_{j,0}(\rho)+\overline{I}_j =0 ,\  \rho<\ell_j,\\
  & D\frac{\partial U_{j,0}(\ell_j^+)}{\partial \rho}=\overline{D}_j\frac{\partial V_{j,0}(\ell_j^-)}{\partial \rho}=\kappa_j [U_{j,0}(\ell_j^+) -V_{j,0}(\ell_j^-)].\end{align}
  \end{subequations}
Assuming that $U_{j,0}(\rho)$ vanishes at infinity, we obtain a solution of the form
\begin{equation}
\label{1Dsola}
U_{j,0}(\rho)=\frac{A_j}{\rho},\quad V_{j,0}(\rho)= \frac{\overline{\calI}_j}{\overline{\gamma}_j}+  B_j \frac{ \sinh(\overline{\beta}_j\rho)}{\rho}, 
\end{equation}
with $\overline{\beta}_j=\sqrt{\overline{\gamma}_j/\overline{D}_j}$. Substituting into the boundary conditions (\ref{sph}c) gives
\begin{align}
-\frac{DA_{j}}{\ell_j}&=\kappa_j \left [A_{j}  -\frac{\overline{\calI}_j\ell_j}{\overline{\gamma}_j} - B_{j} \sinh(\overline{\beta}_j \ell_j) \right ]=\overline{D}_jB_{j}\left \{\overline{\beta}_j\cosh(\overline{\beta}_j\ell_j) - \frac{1}{\ell_j}\sinh(\overline{\beta}_j\ell_j) \right \},
\end{align}
which can be rearranged to yield the result
\begin{equation}
\label{permFj}
A_j=\frac{\overline{D}_j\left \{\overline{\beta}_j\ell_j \cosh(\overline{\beta}_j\ell_j) -  \sinh(\overline{\beta}_j\ell_j) \right \}}{D\ell_j^{-1} \sinh(\overline{\beta}_j\ell_j) +\overline{D}_j\Lambda_j^{-1} \left \{\overline{\beta}_j\ell_j\cosh(\overline{\beta}_j\ell_j)  -\sinh(\overline{\beta}_j\ell_j) \right \}}\frac{\overline{\calI}_j}{\overline{\gamma}_j},
\end{equation}
with $\Lambda_j$ defined in (\ref{Lam}).
Note that we recover equation (\ref{Lam}) with $c_{j,0}=\overline{\calI}_j/\overline{\gamma}_j$ in the fast diffusion limit $\overline{D}_j\rightarrow \infty$. 

\paragraph{\underline{Biochemical reaction network (model III).}}

The 3D version of equation (\ref{stmass}) is
\begin{align}
\label{stmass3D}
\widehat{f}_a(\w_j^*)-  \kappa_j \delta_{a,1} 4\pi \ell_j^2(w_{j,0}^*- U_j(\ell_j)) =0.
\end{align}
Using equation (\ref{Ujr}) with $c_{j,0}=w_{j,0}^*$, we thus obtain the nonlinear algebraic system 
\begin{equation}
\label{Fj3D}
\widehat{f}_a(\w_j^*)+\frac{ 4\pi D\kappa_j\ell_j^2w_{j,0}^*}{ \kappa_j\ell_j+D }\delta_{a,1}+O(\epsilon)=0.
\end{equation}
The existence of steady-state solutions and the reduction of the PDE-ODE system to a nonlinear ODE system is developed in Ref. \cite{Muller13}. However, as far as we are aware, the emergence of collective oscillations via a Hopf bifurcation has not yet been explored in 3D.

\section{Relaxation to steady state and the accumulation time}

Models I and II are linear diffusion problems and we expect the system to relax to a unique steady-state. In many applications one is also interested in characterizing the relaxation to steady state. A standard approach is to calculate the principal nonzero eigenvalue of the negative Laplacian \cite{Ward93,Ward93a,Cheviakov11}. However, this is a global measure of the relaxation process that does not account for possible differences in the relaxation rate at different spatial locations. In addition, it relies on the assumption that the eigenvalues have sufficiently large spectral gaps. Therefore, in this section we use asymptotic methods to calculate the accumulation time for singularly perturbed diffusion problems \cite{Bressloff22,Bressloff23} .

Following Ref. \cite{Bressloff22}, consider the simple motivating example of diffusion along the finite interval, $x\in [0,L]$, with a constant flux $J_0$ at the end $x=0$ and a reflecting boundary at $x=L$. The concentration $u(x,t)$ satisfies the diffusion equation\footnote{In the particular application to morphogenesis, $u(x,t)$ represents the extracellular morphogen concentration gradient along the body axis of a developing embryo, while the boundary flux at $x=0$ is generated by local protein synthesis \cite{Berez10,Berez11,Berez11a,Gordon11}. The constant $\gamma_0$ is an effective degradation rate due to binding of morphogen to cell surface receptors. The spatially varying morphogen concentration drives a corresponding spatial variation in gene expression through some form of concentration thresholding mechanism.This converts a continuously varying morphogen concentration into a discrete spatial pattern of differentiated gene expression across a cell population. One constraint of the model is that the concentration gradient should be established within appropriate developmental time scales. } 
\begin{equation}
\label{grad0}
\frac{\partial u}{\partial t}=-{\mathbb L} u\equiv D\frac{\partial^2u}{\partial x^2}-\gamma_0 u ,\ 0<x<L;\quad \left . -D\frac{\partial u}{\partial x}\right |_{x=0}=J_0, \quad \left . -D\frac{\partial u}{\partial x}\right |_{x=L}=0,
\end{equation}
with $u(x,0)=0$. If $L \gg \xi\equiv \sqrt{{D}/{\gamma_0}}$, then the boundary condition at $x=L$ can be neglected and $u(x)=(J\xi/D)\e^{-x/\xi}$ with $\xi$ identified as the space constant for exponential decay.
One way to quantify the relaxation to steady state is to consider the eigenfunction expansion $0<\lambda_0< \lambda_1\ldots$ 
\begin{equation}
\label{eig}
u(x,t)-u(x)=\sum_{n\geq 0}c_n \phi_n(x)\e^{-\lambda_n t},
\end{equation}
where $(\lambda_n,\phi_n(x))$, $n\geq 0$, are the eigenvalues and eigenfunctions of the linear operator ${\mathbb L}=-D\partial^2_x+\gamma_0$ with homogeneous boundary conditions:
\begin{equation}
{\mathbb L}\phi_n(x) =\lambda_n\phi_n(x),\ x\in [0,L],\quad\phi_n'(0)=\phi_n'(L)=0.
\end{equation}
The eigenfunctions form a complete orthonormal set
\begin{equation}
 \int_0^L\phi_m(x)\phi_n(x)dx=\delta_{n,m}.
 \end{equation}
A straightforward calculation shows that
\begin{equation}
\phi_n(x)=A_n \cos(n\pi (L-x)/L),\quad \lambda_n = \gamma_0+\frac{n^2\pi^2 D}{L^2}, \quad n\geq 0.
\end{equation}
If the positive eigenvalues are well separated, then the relaxation to steady state will be dominated by the term $c_0\phi_0(x)\e^{-\lambda_0 t}$ with $\lambda_0=\gamma_0$ and $\phi_0(x)=$ constant. (The constants $c_n$ are determined by the initial condition.) Hence, we can identify $1/\gamma_0$ as the effective relaxation time.

In order to define the corresponding accumulation time for the above 1D problem, let 
\begin{equation}
\label{accu}
Z(x,t)=1-\frac{u(x,t)}{u(x)}
\end{equation}
denote the fractional deviation of the concentration from steady state. Assuming that there is no overshooting, i. e. no sign reversal of $Z(x,t)$), then $1-Z(x,t)$ can be interpreted as the fraction of the steady-state concentration that has accumulated at $x$ by time $t$. It follows that $-\partial_t Z(x,t)dt$ is the fraction accumulated in the interval $[t,t+dt]$. The accumulation time $T(x)$ at position $x$ is then defined as \cite{Berez10,Berez11,Gordon11}:
\begin{equation}
\label{accu2}
T(x)=\int_0^{\infty} t\left (-\frac{\partial Z(x,t)}{\partial t}\right )dt=\int_0^{\infty} Z(x,t)dt.
\end{equation}
In terms of the eigenfunction expansion (\ref{eig}), 
\begin{equation}
T(x)=-\frac{1}{u(x)}\int_0^{\infty} \sum_{n\geq 0}c_n \phi_n(x)\e^{-\lambda_n t}dt=-\sum_{n=0}^{\infty}\frac{c_n \phi_n(x)}{\lambda_n u(x)},
\end{equation}
which is non-singular since $\lambda_n >0$ for all $n\geq 0$. For the simple diffusion process given by equation (\ref{grad0}), $T(x)$ can be calculated explicitly in the limit $L\rightarrow \infty$ and one finds that \cite{Berez10}
\begin{equation}
\label{morpht}
T(x)=\frac{1}{2\gamma_0}\left (1+\sqrt{\frac{\gamma_0}{D}}x\right ).
\end{equation}

Let us now turn to equations (\ref{gmaster}) with $c_j(\x,t) =c_{j,0}$ (model I) and an initial condition $u(\x,0)=u_0(\x)$ with $u_0(\x)$ having compact support on $\Omega\backslash \calU$. Suppose that there exists a unique solution $u(\x)$ of the steady-state equations (\ref{masterI}). In contrast to the above 1D example, care has to be taken regarding the sign of $u(\x)-u(\x,t)$. If the initial condition $u_0(\x)=0$, then  $u(\x)-u(\x,t)>0$ for all $\x\in \Omega\backslash \calU$ since the boundaries $\calU_j$ act as source terms. On the other hand, if $u_0(\x)$ has sufficient nonzero compact support then the boundaries act as sinks and $u(\x)-u(\x,t)<0$ for all $\x\in \Omega\backslash \calU$. We will consider the latter example here and thus reverse the sign in the 
  definition of the accumulation time:
 \begin{equation}
 T(\x)=\int_0^{\infty} t\left (\frac{\partial Z(\x,t)}{\partial t}\right )dt=-\int_0^{\infty} Z(\x,t)dt=-\lim_{s\rightarrow 0} s\widetilde{Z}(\x,s),
 \end{equation}
where we have introduced the Laplace transform $\widetilde{Z}(\x,s)=\int_0^{\infty}\e^{-st}Z(\x,t)dt$. Laplace transforming the higher-dimensional version of equation (\ref{acc}) implies that
\begin{equation}
s\widetilde{Z}(\x,s)=1-\frac{s\widetilde{u}(\x,s)}{u(\x)},
\end{equation}
so that the corresponding accumulation time can be written as
\begin{eqnarray}
 T(\x)=- \lim_{s\rightarrow 0}\frac{1}{s}\left [1-\frac{s\widetilde{u}(\x,s)}{u(\x)}\right ] =\frac{1}{u(\x)}
\left .\frac{d}{ds}[s\widetilde{u}(\x,s)\right |_{s=0}.
\label{acc}
\end{eqnarray}
(In order to calculate the accumulation time for higher-dimensional diffusion in a singularly perturbed domain, it is more convenient to work in Laplace space.) Laplace transforming equations (\ref{gmaster}) with $c_j(\x,s)=c_{j,0}/s$ and $\calI_0=0$, we have
\begin{subequations} 
\label{LTgmaster}
\begin{align}
D{\bm \nabla}^2\u(\x,s)-(s+\gamma_0)\u(\x,s) +u_0(\x)&=0,\quad \x\in \Omega \backslash   \calU,\\
	D{\bm \nabla} \u(\x,s) \cdot \n(\x)&=0,\  \x\in\partial \Omega ,\\
	D{\bm \nabla} \u(\x,s)\cdot  \n_j(\x) &=\kappa_j \left [\u(\x,s) -\frac{c_{j,0}}{s}\right ],\quad \x\in \partial \calU_j.	\end{align}
	\end{subequations}
Mathematically speaking, these equations are identical in form to the steady-state problem under the mappings $\gamma_0\rightarrow \gamma_0+s$, $c_{j,0}\rightarrow c_{j,0}/s$ and an effective external input $u_0(\x)$. 
Hence, the asymptotic methods presented in section 3 can be extended to solve the diffusion problem in Laplace space. (This also holds for model II, but breaks down in the case of model III due to the presence of nonlinear time-dependent terms; extensions of the latter beyond steady-state analysis will be discussed in section 5.)  Also note that the steady-state solution could also be recovered using the identity
\[u(\x)=\lim_{t\rightarrow \infty} u(\x,t)=\lim_{s\rightarrow 0}s\widetilde{u}(\x,s),\]

One subtle feature about working in Laplace space concerns the small-$s$ behaviour of the corresponding Green's function $G(\x,\x';s)$, which is defined by equations (\ref{GMH}) under the mapping $\gamma_0\rightarrow \gamma_0+s$. (The small-$s$ behaviour is needed in order to evaluate the accumulation time using equation (\ref{acc}).) For $\gamma_0>0$, the Green's function is a nonsingular function of $s$ such that $\lim_{s\rightarrow 0}G(\x,\x';s)=G(\x,\x')$. However, in the special case $\gamma_0=0$ (zero degradation), $G$ has a simple pole at $s=0$ such that \cite{Bressloff22,Bressloff23}
\begin{equation}
\label{Gs}
G(\x,\x';s)=\frac{1}{s|\Omega|}+G_0(\x,\x')+s\partial_s G(\x,\x';0)+O(s^2),
\end{equation}
where $G_0$ satisfies equations (\ref{exG1}). This causes particular problems in 3D, since we obtain a double series expansion in $\epsilon$ and $s$ that includes terms in powers of $\epsilon/s$. Although such terms can be eliminated by performing partial summations, the details are rather technical \cite{Bressloff23}. Therefore, we will restrict the analysis to the case $\gamma_0>0$.

\subsection{Asymptotic analysis in 2D} Following along identical lines to section 3.1, the inner solution takes the form
\begin{align}
\label{LTUV}
    &\widetilde{U}_j(\x,s) =  \frac{c_{j,0}}{s} +\frac{\nu \widetilde{A}_j(s)D}{\kappa_j\ell_j}+\nu \widetilde{A}_j(s)\log[\rho/\ell_j ],\quad \ell_j \leq \rho < \infty ,\quad 
\end{align}
provided that $s<1/\epsilon^2$. (Since we are ultimately interested in the limit $s\rightarrow 0$, this restriction is not an issue here.) Similarly, the outer solution becomes
\begin{align}
    \u(\x,s) &\sim  \Gamma_0(\x,s) - 2\pi \nu D\sum_{k=1}^N \widetilde{A}_k(s)G(\x,\x_k;s),
    \label{LToutu}
\end{align}
where 
\begin{equation}
\Gamma_0(\x,s)=\int_{\Omega} G(\x,\x';s)u_0(\x')dx'.
\end{equation}
Matching the inner and outer solutions yields a system of $N$ equations for the $N$ unknown coefficients $\widetilde{A}_j(s)$:
\begin{align}
  & \widetilde{A}_j (s)\left[1 +  \nu \Psi_j-\nu \log \ell_j + 2\pi  D \nu R(\x_j,\x_j;s )\right]+ 2\pi  D \nu \sum_{k \neq j}\widetilde{A}_k(s) G(\x_j, \x_k;s) \nonumber \\
   &=-\frac{c_{j,0}}{s}+\Gamma_0(\x,s).
  \label{matchacc}
\end{align}
 We thus obtain the matrix solution
 \begin{align}
\label{najacc}
    \widetilde{A}_j(s) 
    &= -\sum_{i=1}^N\left({\bf I}  + \nu {\bf M}(s) \right)^{-1}_{ji}\left [\frac{c_{i,0}}{s} -\Gamma_0(\x_i,s)\right ] ,\quad j = 1, \ldots, N,
    \end{align}
 with
  \begin{equation}
  \label{Ms}
 M_{ji}(s)=\Psi_j\delta_{i,j} +2\pi  D \calG_{ji}(s),
    \end{equation}
    and 
    \begin{equation}
\calG_{ij}(s)=G(\x_i,\x_j;s),\ i\neq j,\quad \calG_{jj}=R(\x_j,\x_j;s)-\frac{\log \ell_j}{2\pi D}.
\end{equation}
 
 For $\gamma_0>0$, the Green's function $G(\x,\x';s)$ is non-singular in the limit $s\rightarrow 0$. Hence, multiplying the outer solution (\ref{LToutu}) by $s$, differentiating with respect to $s$ and taking the limit $s\rightarrow 0$ gives
\begin{align}
   \lim_{s\rightarrow 0} \frac{ds \u(\x,s)}{ds} &\sim   - 2\pi \nu D \sum_{k=1}^N\left \{ A_k'G(\x,\x_k)+A_kH(\x,\x_k)\right \}+\Gamma_0(\x,0),
\end{align}
where 
\begin{equation}
H(\x,\x_k)=\lim_{s\rightarrow 0}\frac{\partial G(\x,\x_k;s)}{\partial s},
\end{equation}
$A_j=\lim_{s\rightarrow 0}s\widetilde{A}_j(s)$ is the coefficient of the corresponding steady-state solution, and
\begin{equation}
A_j'=\lim_{s\rightarrow 0}\left ( \frac{d[s\widetilde{A}_j(s)]}{ds} \right )=\nu \sum_{i=1}^N  {\bm \Theta}_{ji}c_{i,0} +\sum_{i=1}^N\left({\bf I}  + \nu {\bf M}(0) \right)^{-1}_{ji} \Gamma_0(\x_i,0),
\end{equation}
with
\begin{equation}
{\bm \Theta}=\left({\bf I}  + \nu {\bf M}(0) \right)^{-1} {
\bm M}'(0)\left({\bf I}  + \nu {\bf M}(0) \right)^{-1}.
\end{equation}
We have multiplied both sides of equation (\ref{najacc}) by $s$ and differentiated with respect to $s$ using the matrix identity
\begin{equation}
0=\frac{d}{ds}\left [{\bf T}(s) {\bf T}(s)^{-1}\right ] ={\bf T}'(s) {\bf T}(s)^{-1}+{\bf T}(s)\frac{d}{ds}{\bf T}(s)^{-1}
\end{equation}
for any invertible and differentiable matrix ${\bf T}(s)$. Combining the various results yields the following general expression for the accumulation time:
\begin{eqnarray}
 T(\x)&\sim&\frac{\Gamma_0(\x,0)-2\pi \nu D \sum_{k=1}^N\left \{ A_k'G(\x,\x_k)+A_kH(\x,\x_k)\right \}}{2\pi \nu D \sum_{k=1}^NA_kG(\x,\x_k)}\nonumber \\
 &\sim&\frac{\int_{\Omega}G(\x,\x')u_0(\x')d\x'}{2\pi \nu D \sum_{k=1}^NG(\x,\x_k)c_{k,0}}+O(1).
 \label{acc2D}
\end{eqnarray}

A general result of this analysis is that the accumulation time diverges as $\nu \rightarrow 0$ in the case of a nonzero initial condition $u_0(\x)$ in the bulk, for which the compartments act as strongly localised sinks. The situation is very different if $u_0(\x)\equiv 0$ for all $\x\in \Omega\backslash \calU$, since the compartments now act as strongly localised sources from which particles diffuse into the bulk. Reversing the sign in the definition of the accumulation time, we have $T(\x)=O(1)$.
\subsection{Asymptotic analysis in 3D}

We now turn to the aysmptotic analysis of equations (\ref{LTgmaster}) in 3D. 
The outer solution can be expanded as
\begin{equation}
\u(\x,s)\sim \Gamma_0(\x,s)+ \epsilon \u_1(\x,s)+\epsilon^2 \u_2(\x,s)+\ldots
\end{equation}
with
\begin{subequations}
\label{LT3Douter}
\begin{align}
&D{\bm \nabla}^2 \u_{n}(\x,s) -\gamma_0 u_{n}(\x,s)=0,\quad \x\in \Omega\backslash \{\x_1,\ldots,\x_N\},\\
& {\bm \nabla}\u_{n}(\x,s)\cdot \n(\x)=0,\quad \x \in \partial \Omega
\end{align}
\end{subequations}
for $n=1,2.\ldots$. Equations (\ref{LT3Douter}) are supplemented
by singularity conditions as $\x\rightarrow \x_j$, $j=1,\ldots,N$. The latter are determined by matching to the inner solution
\begin{equation}
\widetilde{U}_j(\y,s) = \widetilde{U}_{j,0}(\y,s) +\epsilon \widetilde{U}_{j,1}(\y,s)+\ldots,
\end{equation}
where
\begin{subequations}
\label{LT3Dinner}
\begin{align}
&D {\bm \nabla}^2_{\y}\widetilde{U}_{j,n}(\y,s)= \gamma_0 \widetilde{U}_{j.n-2}(\y,s),\quad |\y| > \ell_j,\\
	&D{\bm \nabla}_{\y} \widetilde{U}_{j,n}(\y,s)\cdot \n_j - \kappa_j \widetilde{U}_{j,n}(\y,s)=-\frac{\kappa_j c_{j,0}}{s}\delta_{n,0},\quad   |\y |=\ell_j  .
\end{align}
\end{subequations}
The matching condition is now
\begin{equation}
\label{LT3Dmatch}
 \Gamma_0(\x,s)+ \epsilon \u_1(\x,s)+\epsilon^2 \u_2(\x,s) \sim \widetilde{U}_{j,0}(\x,s)+  \epsilon \widetilde{U}_{j,1}(\x,s)+\ldots
\end{equation}
 In order to perform this matching, we need 
to Taylor expand $\Gamma(\x,s)$ near
the $j$-th compartment and rewrite it in terms of
stretched coordinates:
\begin{equation}
\label{p0}
\Gamma_0(\x,s)  \sim \Gamma_0(\x_j,s)+ \epsilon {\bm \nabla}\Gamma_0(\x_j,s) \cdot \y +\ldots
\end{equation}
It immediately follows that the $O(1)$ contribution to the inner solution no longer vanishes at infinity. In particular,
\begin{equation}
\widetilde{U}_{j,1}(\y,s)=\Gamma_0(\x_j,s)+\frac{\Lambda_j }{\rho}  \left [ \frac{c_{j,0}}{s}-\Gamma_0(\x_j,s)\right ].
\end{equation}
The matching condition (\ref{LT3Dmatch}) then implies that $u_1(\x)$ satisfies equations (\ref{LT3Douter}) together with the singularity condition
\[\u_1(\x,s)\sim \frac{\Lambda_j }{|\x-\x_j|} \left [ \frac{c_{j,0}}{s}-\Gamma_0(\x_j,s)\right ]\quad \mbox{as } \x\rightarrow \x_j.\]
Hence, we have the solution
\begin{align}
\u_1(\x,s) 
 \sim {4\pi}  \epsilon D \sum_{k=1}^N \Lambda_k\left [ \frac{c_{k,0}}{s}-\Gamma_0(\x_k,s)\right ] G(\x,\x_k;s) .
\end{align}
Combining the non-singular near-field behaviour of $\u_1(\x,s) $ with equation (\ref{p0}) implies that
\begin{equation}
\widetilde{U}_{j,2}(\y,s)\rightarrow {\bm \nabla}\Gamma_0(\x_j,s) \cdot \y+{4\pi}   D \sum_{k=1}^N \Lambda_k\left [ \frac{c_{k,0}}{s}-\Gamma_0(\x_k,s)\right ] \calG_{jk}(s)\mbox{  as. } |\y|\rightarrow \infty,
\end{equation}
where
${\mathcal G}_{ij}(s) =G(\x_i,\x_j;s$ for $i\neq j$, and ${\mathcal G}_{ii}(s) =R(\x_i,\x_i;s)$.

The first term on the right-hand side generates contributions to the inner solution in the form of first-order spherical harmonics, which only affect the outer solution at $O(\epsilon^3)$. In order to show this, we
decompose the solution around the $j$-th compartment as $U_{j,1}=A_{j,1}+B_{j,1}$ with 
\[A_{j,1}(\y,s)\rightarrow \chi_j (s)={4\pi}   D \sum_{k=1}^N \Lambda_k\left [ \frac{c_{k,0}}{s}-\Gamma_0(\x_k,s)\right ] \calG_{jk}(s) \mbox{ as } |\y|\rightarrow \infty,\]
and
\[B_{j,1}(\y,s)\rightarrow {\bf b}_j\cdot \y \mbox{ as } |\y|\rightarrow \infty,\quad {\bf b}_j(s)={\bm \nabla}\Gamma_0(\x_j,s) .\]
The solution for $A_{j,1}$ is then 
\begin{equation}
A_{j,1}(\y,s)=\chi_j(s)\left (1-\frac{\Lambda_j}{|\y|}\right ).
\end{equation}
In order to determine $B_{j,1}$, we introduce local spherical polar coordinates such that ${\bf b}_j=(0,0,b_j)$ and $\y\cdot {\bf b}_j=b_jr\cos \theta$, $0\leq \theta \leq \pi$. In spherical polar coordinates we have
\begin{align}
&\frac{\partial^2B_{j,1}}{\partial r^2}+\frac{2}{r}\frac{\partial B_{j,1}}{\partial r} +\frac{1}{r^2\sin \theta}\frac{\partial }{\partial \theta}\left (\sin \theta\frac{\partial B_j^{(1)}}{\partial \theta} \right )=0,\  r>1,\\
&   B_{j,1} \sim b_jr\cos \theta  \mbox{ as } r \to \infty;\quad B_{j,1}=0 \mbox{ on } r=\ell_j. 
\end{align}
Recall that Laplace's equation in spherical polar coordinates has the general solution
\begin{equation}
\label{gen}
B(r,\theta,\phi) = \sum_{l\geq 0}\sum_{m=-l}^l \left (a_{lm}r^l+\frac{b_{lm}}{r^{l+1}}\right )P_l^m(\cos \theta)\e^{im\phi},
\end{equation}
where $P_l^m(\cos \theta)$ is an associated Legendre polynomial. Imposing the Dirichlet boundary condition and the far-field condition implies that
\begin{equation}
\label{B1}
B_{j,1}(\y,s)=b_j\ell_j\cos \theta \left (\frac{|\y|}{\ell_j}-\frac{\ell_j^2}{|\y|^2}\right ).
\end{equation}
Clearly this leads to a singular term that is $O(\epsilon^3)$.

In summary, we obtain an outer solution of the form
\begin{align}
\u(\x,s) 
 \sim \Gamma_0(\x,s)+{4\pi}  \epsilon D \sum_{k=1}^N \Lambda_k\left [ \frac{c_{k,0}}{s}-\Gamma_0(\x_k,s)-\epsilon \chi_j(s)\right ]  G(\x,\x_k;s) .
\end{align}
Multiplying both sides by $s$, differentiating with respect to $s$ and taking the limit $s\rightarrow 0$ gives
\begin{align}
   \lim_{s\rightarrow 0} \frac{ds \u(\x,s)}{ds} &\sim    4\pi \epsilon D \sum_{k=1}^N \Lambda_k
  \left \{  c_{k,0}H(\x,\x_k)-\Gamma_0(\x_k,0)G(\x,\x_k)\right \}+\Gamma_0(\x,0),
    \label{LToutu2}
\end{align}
We conclude that
\begin{eqnarray}
 T(\x)&\sim&\frac{\Gamma_0(\x,0)+ 4\pi \epsilon D \sum_{k=1}^N \Lambda_k
  \left \{  c_{k,0}H(\x,\x_k)-\Gamma_0(\x_k,0)G(\x,\x_k)\right \}}{4\pi  D \epsilon \sum_{k=1}^N\Lambda_kc_{k,0}G(\x,\x_k) }\nonumber \\
 &\sim&\frac{\int_{\Omega}G(\x,\x')u_0(\x')d\x'}{4\pi  D \epsilon \sum_{k=1}^N\Lambda_kc_{k,0}G(\x,\x_k) }+O(1).
\end{eqnarray}
As in the 2D case, $T(\x)$ is singular in the limit $\epsilon \rightarrow 0$ for $u_0(\x)$ having nonzero support, but is finite when $u_0\equiv 0$ (after a sign reversal).

\section{Asymptotic reduction of a 2D PDE-ODE model to a nonlinear ODE system}

In contrast to the linear models I and II , which are based on linear diffusion equations, model III describes a nonlinear coupled PDE-ODE system that can support multiple steady-state solutions. However, only those that are linearly stable will be observable. Moreover, a steady-state solution may undergo a bifurcation due to one or more eigenvalues associated with the linearised system crossing the imaginary axis, resulting in a modified set of unstable and stable steady-state solutions, and possibly the emergence of a limit cycle. The linear stability of the 2D model and conditions for a Hopf bifurcation have been analysed in Ref. \cite{Gou16}. Asymptotic methods have also been used to reduce the PDE-ODE system to an effective nonlinear ODE system in the well mixed limit $D\gg O(\nu^{-1})$ \cite{Gou16,Iyan21a}, and this analysis has subsequently been extended to the regime $D=O(\nu^{-1})$ \cite{Iyan21b}. The reduced models provided a more tractable mathematical framework for investigating the role of diffusion-mediated coupling on the collective dynamics of compartmental oscillators, at least in the fast diffusion regime. 
In this section we review the asymptotic reduction of Ref. \cite{Iyan21b}\footnote{Extending asymptotic and numerical methods to the time domain is an important topic for future work in this area. A recent preprint by Pelz and Ward \cite{Pelz24} develops a hybrid asymptotic-numerical theory for studying diffusion-mediated collective oscillations. In particular, these authors extend asymptotic methods for analyzing equations (\ref{gmaster}) and (\ref{mass}) to the time-dependent setting without any restrictions on the diffusivity. This leads to a novel integro-differential ODE system that characterizes intracellular dynamics in a memory-dependent bulk-diffusion field. They also introduce a fast time-marching scheme to compute numerical solutions to the integro-differential system over long time intervals.}.

The main idea is to derive an ODE for the mean bulk concentration
\begin{equation}
\overline{u}(t)=\frac{1}{|\Omega\backslash\calU|}\int_{\Omega\backslash \calU}u(\x,t)d\x,
\end{equation}
which is accurate to $O(\nu)$ and couples to the dynamical variables $\w_j(t)$, $j=1,\ldots,N$. First, integrating equation (\ref{gmaster}a) with $\kappa_j\rightarrow \kappa_j/\epsilon$ and $\calI_0=0$, and using the divergence theorem shows that
\begin{subequations}
\label{wu}
\begin{equation}
\frac{d \overline{u}(t)}{\partial t} +\gamma_0 \overline{u}(\x,t)= \frac{2\pi}{|\Omega|}\sum_{j=1}^N \ell_j\kappa_j\left [w_{j,0}(t)-\frac{1}{2\pi \epsilon \ell_j}\int_{\partial \calU_j} u(\x,t)d\x \right ],
\end{equation}
where we have taken $|\Omega\backslash\calU|=|\Omega|+o(\epsilon)$ (assuming that $N$ is not too large).
Moreover,
\begin{align}
 |\calU_j|\frac{dw_{j,a}}{dt}=\widehat{f}_a(\w_j)- 2\pi  \kappa_j \ell_j\delta_{a,1}\left [w_{j,0}(t)-\frac{1}{2\pi \epsilon \ell_j}\int_{\partial \calU_j} u(\x,t)d\x \right ].
\end{align}
\end{subequations}
Assuming $D=O(\nu^{-1})$, set $D=D_0/\nu$ with $D_0=O(1)$. The inner solution around the $j$-th compartment satisfies the equation
\begin{subequations}
\begin{align}
&D {\bm \nabla}^2_{\y}U_j(\y,t)= 0,\quad |\y| > \ell_j,\\
	&{\bm \nabla}_{\y} U_j(\y,t)\cdot \n_j(\y)= \frac{\nu \kappa_j }{D_0}[U_j(\y,t)-w_{j,0}(t)],\quad   |\y |=\ell_j  .
\end{align}
\end{subequations}
Using polar coordinates with $|\y|=\rho$ and setting $U_j=U_j(\rho,t)$, we have the solution
\begin{align}
\label{UVQS}
    &U_j(\rho,t) =  \Phi_j(t)+\nu A_j(t) \log[\rho/\ell_j ],\quad \ell_j \leq \rho < \infty ,
    \end{align}
    where
    \begin{equation}
A_j(t)=\frac{\kappa_j\ell_j}{D_0}[\Phi_j(t)-w_{j,0}(t)].
\end{equation}
Substituting the inner solution into equations (\ref{wu}) then gives
\begin{subequations}
\label{wu2}
\begin{align}
&\frac{d \overline{u}(t)}{\partial t} +\gamma_0 \overline{u}(\x,t)=- \frac{2\pi D_0}{|\Omega|}\sum_{j=1}^NA_j(t),\\
& |\calU_j|\frac{dw_{j,a}}{dt}-\widehat{f}_a(\w_j) =2\pi  D_0 \delta_{a,1}A_j(t).
\end{align}
\end{subequations}

It remains to determine the coefficients $A_j(t)$, $j=1,\ldots,N$ by matching the inner solution with the outer solution, which satisfies
\begin{subequations} 
\label{outerQS}
\begin{align}
	\frac{\partial u(\x,t)}{\partial t} &= \frac{D_0}{\nu}{\bm \nabla}^2u(\x,t)-\gamma_0 u(\x,t) ,\quad \x\in \Omega \backslash \{\x_1,\ldots,\x_N\},\\
	D{\bm \nabla} u(\x,t) \cdot \n(\x)&=0,\  \x\in\partial \Omega ,\\
	u(\x,t)&\sim w_{j,0}(t)+  \nu A_j(t)  \log(|{\bf x} - {\bf x}_j|/\ell_j) +\left [1+\frac{D_0}{\kappa_k\ell_j}\right]A_j(t).
		\end{align}
	\end{subequations}
Introduce the expansion
\begin{equation}
u(\x,t)\sim \overline{u}(t)+\frac{\nu}{D_0}u_1(\x,t)+\ldots,\quad \mbox{where  }\int_{\Omega}u_1(\x,t)d\x=0. 
\end{equation}
 Substituting into the outer equations shows that \cite{Iyan21b}
 \begin{subequations} 
\label{outerQS1}
\begin{align}
	& D_0{\bm \nabla}^2u_1(\x,t)=\frac{d\overline{u}(t)}{dt}+\gamma_0 \overline{u}(t)+2\pi \sum_{k=1}^NA_k(t)\delta(\x-\x_k)  ,\quad \x\in \Omega ,\\
&	D{\bm \nabla} u(\x,t) \cdot \n(\x)=0,\  \x\in\partial \Omega ,\\
&	u_1(\x,t)\sim \frac{w_{j,0}(t)-\overline{u}(t)}{\nu}+    A_j(t)  \log(|{\bf x} - {\bf x}_j|/\ell_j) +\left [1+\frac{D_0}{\kappa_j\ell_j}\right]\frac{A_j(t)}{\nu}.
		\end{align}
	\end{subequations}
	The sum over Dirac delta functions in equation (\ref{outerQS1}) ensures that is consistent with equation (\ref{wu2}) when integrated over $\Omega$. Proceeding as in previous sections, we write
\begin{equation}
u_1(\x,t)\sim -2\pi D_0\sum_{k=1}^NA_k(t)G_0(\x,\x_k),
\end{equation}
where $G_0$ satisfies equations (\ref{exG1}). Finally, matching the inner and outer solutions using equation (\ref{outerQS1}c)  gives
\begin{align}
   \left [1+\frac{D_0}{\kappa_j\ell_j}\right]A_j (t)+ 2\pi  D_0 \nu \sum_{k=1}^N \calG_{0,jk}A_k(t) =\overline{u}(t)-w_{j,0}(t),\quad j=1,\ldots,N.
  \label{matchQS}
\end{align}

In summary, as originally shown in Ref. \cite{Iyan21b}, the PDE-ODE system of model III with bulk diffusivity $D=D_0/\nu =O(\nu^{-1})$ can be reduced to an ODE system given by equations (\ref{wu2}) with coefficients $A_j(t)$ determined by the matrix equation (\ref{matchQS}). In order to relate this result to the classical ODE system (\ref{mass0}), suppose that all the cells are identical so that $\kappa_j=\kappa$ and $\ell_j=\ell$. The solution of equation (\ref{matchQS}) can then be written as \begin{equation}
A_j(t)=\frac{\kappa \ell}{\kappa\ell+D_0}\sum_{k=1}^N[{\bf I}+\nu {\bf Q}]^{-1} _{jk} [\overline{u}(t)-w_{k,0}(t)]\quad \mbox{where } {\bf Q}=\frac{\kappa \ell D_0}{\kappa \ell+D_0}{\bm \calG}_0,
\end{equation}
and equations (\ref{wu2}) become
\begin{subequations}
\label{wu3}
\begin{align}
&\frac{d \overline{u}(t)}{\partial t} +\gamma_0 \overline{u}(\x,t)=\frac{2\pi }{|\Omega|}\frac{\kappa \ell D_0}{\kappa\ell+D_0}\sum_{j,k=1}^NW_{jk} [w_{k,0}(t)-\overline{u}(t)],\\
& \frac{dw_{j,a}}{dt}={f}_a(\w_j) -\frac{ 2\pi N}{|\calU|} \frac{\kappa \ell D_0}{\kappa\ell+D_0} \delta_{a,1}\sum_{j=1}^NW_{jk}[w_{k,0}(t)-\overline{u}(t)],\quad a=1,\ldots,K,
\end{align}
\end{subequations}
where ${\bf W}=[{\bf I}+\nu {\bf Q}]^{-1} $. Since ${\bf W}\rightarrow {\bf I}$ in the well-mixed limit $D_0\rightarrow \infty$, equations (\ref{wu3}) reduce to the non-spatial model (\ref{mass0}) with $\widehat{\kappa}=2\pi N\kappa \ell/|\calU|$. However, all information regarding the scaled diffusivity $D_0$ and the spatial configurations $\x_1,\ldots,\x_N$ of the cells via the dependence of ${\bf W}$ on $G_0(\x_i,\x_j)$ is lost.  The effects of these features on synchrony and oscillatory dynamics of diffusion-mediated quorum sensing are explored in some detail in Ref. \cite{Iyan21b} by considering the particular example of Sel'kov reaction kinetics. (The latter was originally introduced as a simple model
 of glycolysis oscillations in yeast cells \cite{Selk68}.) Each isolated compartment is assumed to be slightly below the threshold for the onset of oscillations. In the presence of diffusion-mediated cell--cell interactions with $D_0=O(1)$, one finds that for a range of parameter values and geometric arrangements of cells, the coupled system can undergo a Hopf bifurcation resulting in collective synchronisation. However, these collective oscillations disappear in the well-mixed limit, indicating the crucial role of bulk diffusion in the emergence of synchrony.

\medskip

\section{Themes and variations}

We conclude by discussing some generalisations and open problems arising from the theoretical framework presented in previous sections.

\subsection{Non-spherical compartments}
As originally shown by Ward and Keller \cite{Ward93,Ward93a}, it is possible to generalise the asymptotic analysis of singularly perturbed diffusion problems to more general compartmental shapes such as ellipsoids by applying classical results from electrostatics. In particular, note that the low-order terms in the expansion of the inner solution around the $j$-th compartment satisfy Laplace's equation. This is supplemented by a Robin boundary condition on $\calU_j$ together with some far-field condition. For the sake of illustration, suppose that the Robin boundary condition is replaced by the Dirichlet condition $U(\y)=\Phi$ on $\partial \calU_j$ and $U(\y)\rightarrow u_{\infty}$ as $|\y|\rightarrow \infty$. (Given the solution to Laplace's equation, one could determine $\Phi$ self-consistently by imposing the Robin boundary condition.) We can then write the solution as
\begin{equation}
U(\y)=\Phi+[u_{\infty}-\Phi](1-w(\y)),
\end{equation}
 with $w(\y)$ satisfying the boundary value problem
\begin{align}
\nabla_{\bf y}^2 w(\y)&=0,\ \y\notin \calU_k; \quad w(\y)=1,\ \y \in \partial \calU_j; \quad 
w(\y)\rightarrow 0\quad \mbox{as } |\y|\rightarrow \infty.
\end{align}
This is a well-known problem in electrostatics and has the far-field behaviour
\begin{equation}
\label{eq6:ws}
{w(\y)\sim \frac{C_j}{|\y|}+\frac{{\bf P}_j\cdot \y}{|\y|^3}+\ldots \mbox{as } |\y|\rightarrow \infty,}
\end{equation}
where $C_j$ is the capacitance and ${\bf P}_j$ the dipole vector of an equivalent charged conductor with the shape $\calU_j$. (Here $C_j$ has the units of length.) Some examples of capacitances for various trap shapes are as follows \cite{Cheviakov11}:
\begin{align*}
C_j&=a \mbox{ (sphere of radius } a),\\
C_j&= 2a(1-1/\sqrt{3}) \mbox{ (hemisphere of radius } a),\\
C_j&= \frac{\sqrt{a^2-b^2}}{\cosh^{-1}(a/b)} \mbox{ (prolate spheroid with semi-major and minor axes } a,b),\\
C_j&= \frac{\sqrt{a^2-b^2}}{\cosh^{-1}(b/a)} \mbox{ (oblate spheroid with semi-major and minor axes } a, b).
\end{align*}
One typically finds that low-order terms in the asymptotic expansion for non-spherical compartments can be obtained from the corresponding expressions for the spherical case by replacing the scaled radius $\ell_j$ with the more general shape capacitance $C_j$. An analogous result holds for 2D with $\ln \ell_j $ replaced by the so-called logarithmic capacitance.

    \subsection{Switching processes}
    
    \begin{figure}[b!]
  \begin{center}
  \includegraphics[width=8cm]{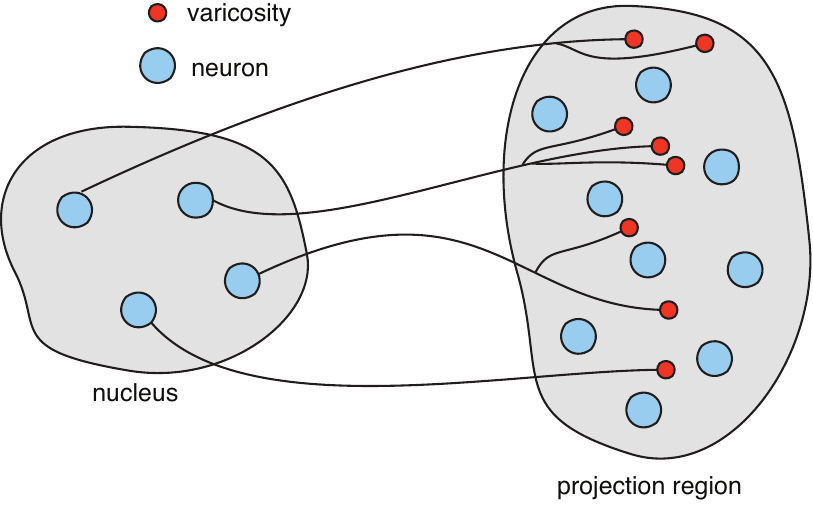}
  \caption{Schematic diagram illustrating volume transmission. Neurons in a brain stem nucleus project their axons to another brain region, where they release neurotransmitter
from numerous small varicosities into the extracellular space.}
  \label{vol}
  \end{center}
\end{figure}

\noindent \underline{Volume transmission.}
There is at least one example of a singularly perturbed diffusion problem in cell biology that is not covered by the boundary conditions shown in Fig. \ref{fig2}, namely, volume transmission. The latter refers to a form of neural communication whereby a population of neurons makes non-specific connections to other brain regions \cite{Fuxe10,Reed23}. That is, rather than forming
one-to-one synaptic connections to other neurons, they simply release
neurotransmitter into the extracellular space from numerous axon terminals known as varicosities.
A number of well-known neurotransmitters participate in volume transmission, including serotonin, dopamine, histamine and acetylcholine. In most cases, the cell bodies are located in brain stem nuclei and their axons project to
a distant brain region such as the striatum or cortex, where they release the neurotransmitter, see Fig. \ref{vol}. 
Recently, volume transmission has been formulated as an example of diffusion in a singularly perturbed domain in which the boundaries $\partial \calU_j$ randomly switch between open and closed states \cite{Lawley18,Lawley20}. That is, each compartment $\calU_j$ represents a varicosity that acts as a source of neurotransmitter when the source neuron fires and is a sink for neurotransmitter otherwise. The latter is due to the re-uptake of neurotransmitter into the terminals.

The model equations take the form
\begin{subequations} 
\label{VT}
\begin{align}
	\frac{\partial u(\x,t)}{\partial t} &= D{\bm \nabla}^2u(\x,t),\quad \x\in \Omega \backslash   \calU,\\
	D{\bm \nabla} u(\x,t) \cdot \n(\x)&=0,\  \x\in\partial \Omega ,\\
	D{\bm \nabla} u(\x,t)\cdot  \n_j(\x) &=\calJ_j, \quad \x\in \partial \calU_j \mbox{ if neuron $j$ is quiescent at time $t$},\\	
	u(\x,t)&=0, \quad \x\in \partial \calU_j \mbox{ if neuron $j$ is firing at time $t$.}\end{align}
	\end{subequations}
	Here $\calJ_j$ is a constant flux of neurotransmitter into the extracellular space when the corresponding source neuron is firing. Note that the volume $\Omega$ is the extracellular space that includes the varicosities but excludes the target cells. Equations (\ref{VT}) are an example of a diffusion process in a randomly switching environment \cite{Bressloff15a,Bressloff15b}. 
One way to analyze such a system is to take expectations with respect to the random switching process, conditioned on the current state of the environment. For the sake of illustration, suppose that all the source neurons fire synchronously and that switching between the global quiescent and firing states is described by a two-state Markov chain with transition rates $\alpha,\beta$:
\begin{equation}
\mbox{(quiescent)  } 0 \Markov{\alpha}{\beta} 1 \mbox{  (firing)}.
\end{equation}
(If each varicosity switched independently then we would have a Markov chain of size $2^N$, where $N$ is the number of varicosities.)
Let $N(t)\in \{0,1\}$ denote the current state of the source neurons, and introduce the conditional expectations
\begin{equation}
\overline{u}_n(\x,t)=\E[u(\x,t){\bm 1}_{N(t)=n}].
\end{equation}
It can then be shown that $u_n$ evolves according to the equations
\begin{subequations} 
\label{VT2}
\begin{align}
	\frac{\partial \overline{u}_n(\x,t)}{\partial t} &= D{\bm \nabla}^2 \overline{u}_n(\x,t)+\sum_{m=01,1}Q_{nm}\overline{u}_m(\x,t),\quad \x\in \Omega \backslash   \calU,\\
	D{\bm \nabla} \overline{u}_n(\x,t) \cdot \n(\x)&=0,\  \x\in\partial \Omega ,\\
	\overline{u}_0(\x,t)&=0,\quad D{\bm \nabla} \overline{u}_1(\x,t)\cdot  \n_j(\x) =\calJ_j, \quad \x\in \partial \calU_j .\end{align}
	\end{subequations}
	The matrix generator of the Markov chain is
	\begin{equation}
	{\bf Q}=\left (\begin{array}{cc} -\beta &\alpha \\ \beta &-\alpha \end{array}  \right ) .
	\end{equation}
The steady-state can then be obtained by performing two successive asymptotic expansions based on model I. First, suppose that we replace the inhomogeneous Neumann boundary condition for $\overline{u}_1$ on $\partial \calU_j$ by the Dirichlet boundary condition $\overline{u}_1(\x,t)=\phi_j$ for some unknown constant $\phi_j$. Summing the steady-state version of equations (\ref{VT2}) with respect to $n=0,1$ and setting $\overline{u}(\x)=\overline{u}_0(\x)+\overline{u}_1(\x)$ gives
 \begin{subequations} 
\label{VT3}
\begin{align}
	D{\bm \nabla}^2\overline{u}(\x)&=0,
	\quad \x\in \Omega \backslash   \calU,\\
	D{\bm \nabla} \overline{u}(\x) \cdot \n(\x)&=0,\  \x\in\partial \Omega, \quad \overline{u}(\x)=\phi_j \quad \x\in \partial \calU_j .\end{align}
	\end{subequations}
	This is equivalent to equations (\ref{masterI}) for model I in the limits $\kappa_j\rightarrow \infty$ and $\gamma_0\rightarrow 0$. Solving the latter equation along the lines outlined in section 3, we then set $n=1$ and $\overline{u}_0(\x)=\overline{u}(\x)-\overline{u}_1(\x)$ in the steady-state version of equations (\ref{VT2}):
	\begin{subequations} 
\label{VT4}
\begin{align}
	& D{\bm \nabla}^2 \overline{u}_1(\x)-(\alpha+\beta)\overline{u}_1(\x)=-\beta \overline{u}(\x),\quad \x\in \Omega \backslash   \calU,\\
	&D{\bm \nabla} \overline{u}_1(\x) \cdot \n(\x)=0,\  \x\in\partial \Omega ,\quad
	\overline{u}_1(\x,t)=\phi_j, \quad \x\in \partial \calU_j .\end{align}
	\end{subequations}
	Equations (\ref{VT4}) yield another version of model I, with $\gamma_0=\alpha+\beta$ and an external input $\beta \overline{u}(\x)$, which can also be solved using the methods of section 3. The resulting solution depends on the unknown constants $\phi_j$, which are finally determined by imposing the original Neumann boundary condition for $\overline{u}_1$ \cite{Lawley20}. (An analogous pair of BVPs arises in the analysis of narrow escape through a stochastically-gated boundary \cite{Bressloff15b}.)
	
One general result that emerges from the singular-perturbation theory is that the leading order term in the asymptotic expansion of the extracellular concentration is a constant that is independent of
the number and location of the varicosities, correlations in neural firings, and the size and geometry of
the extracellular space \cite{Lawley20}. These features do appear in the next leading order term, which may still be space-independent under mild assumptions.	

	\medskip

\noindent \underline{Active phase separation.} Another example of a singularly perturbed diffusion problem with switching is a model for the active suppression of Oswald ripening during the formation of multiple biological condensates. In this model, solute molecules randomly and independently switch between two conformational states, one phase separating ($P$) and the other soluble $(S)$ \cite{Wurtz18,Lee19}. It is assumed that switching between the states $P$ and $S$ occurs according to the chemical reactions
\[P\Markov{h}{k} S,\]
where $h$ and $k$ are concentration-independent reaction rates. The latter reflects the non-equilibrium nature of the chemical reactions, in which detailed balance does not hold due to the phosphorylating action of ATP, say. Let $U_n(\x,t)$, $\x\in \Omega\backslash \calU$, denote the concentration of solute molecules in state $n\in \{S,P\}$. The steady-state diffusion equations take the form \cite{Bressloff20a,Bressloff20b}
\begin{subequations}
\label{outATP}
\begin{align}
D\nabla^2 U_S -k U_S+h U_P&=0, \\
D\nabla^2 U_P+kU_S-hU_P &=0, \quad \x \in \Omega\backslash \calU,
\end{align}
supplemented by the boundary conditions
\begin{align}
{ \bm \nabla}U_S\cdot \n(\x) &=0={ \bm \nabla}U_P\cdot \n(\x) \mbox{ on }   \partial \Omega,\nonumber \\
 U_S(\x)=\Theta_j,\quad U_P(\x)&=c_0(r_j)\equiv \phi_a\left (1+\frac{\ell_c}{r_i}\right ),\,  \mbox{ on } \partial \calU_j.
\end{align}
\end{subequations}
The unknown constants $\Theta_j$ are determined by solving a second pair of diffusion equations within each droplet $\calU_j$ and imposing flux continuity of the soluble phase across $\partial \calU_j$:
\begin{subequations}
\label{inATP}
\begin{align}
D\nabla^2 V_{j,S} -k V_{j,S}+h V_{j,P}&=0, \\
D\nabla^2 V_{j,P}+kV_{j,S}-hV_{j,P} &=0, \quad \x \in \calU_j,
\end{align}
supplemented by the boundary conditions
\begin{align}
 V_{j,S}(\x)=\Theta_j,\quad V_{j,P}(\x)&=\phi_b ,\,  \mbox{ on } \partial \calU_j.
\end{align}
\end{subequations}
The flux continuity conditions are
\begin{equation}
{ \bm \nabla}U_S(\x) \cdot \n(\x) ={ \bm \nabla}V_{j,S}\cdot \n(\x) \mbox{ on }   \partial \calU_j.
\end{equation}

	\begin{figure}[t!]
  \begin{center}
  \includegraphics[width=8cm]{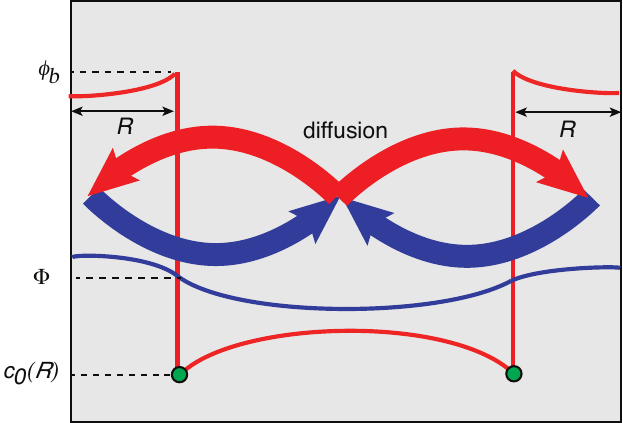}
  \caption{Schematic illustration of the non-equilibrium concentration profiles of $P$-state and $S$-state solute molecules for a pair of liquid droplets in the presence of chemical reactions that drive switching between the two states. Coexistence of the two droplets is maintained by a circulating flux of $P$ and $S$ molecules. Within the $P$-rich droplets, the chemical reaction $P\rightarrow S$ dominates, leading to an accumulation of $S$ molecules inside the droplets, and a net out-flux of $S$ molecules. In the cytoplasm the reverse reaction dominates, leading to the creation and accumulation of $P$ molecules between the droplets and a net in-flux of $P$ molecules. }
  \label{fig8}
  \end{center}
\end{figure}

 The combination of non-equilibrium chemical reactions, phase separation and diffusion can lead to a stationary state consisting of multiple high density droplets that is maintained by circulating fluxes of $P$ and $S$ molecules. The basic mechanism is illustrated in Fig. \ref{fig8} for a pair of droplets. As in classical Ostwald ripening, dynamical equations for the droplet radii can be derived using mean field theory under an adiabatic approximation. The resulting equations support a stable fixed point that represents a multi-droplet solution \cite{Wurtz18,Lee19}. Singular perturbation theory has subsequently been used to investigate corrections to mean field theory,
 which take into account finite-size effects associated with the boundary of the domain and the positions of the droplets \cite{Bressloff20a}. 
\medskip

\begin{figure}[t!]
  \begin{center}
  \includegraphics[width=8cm]{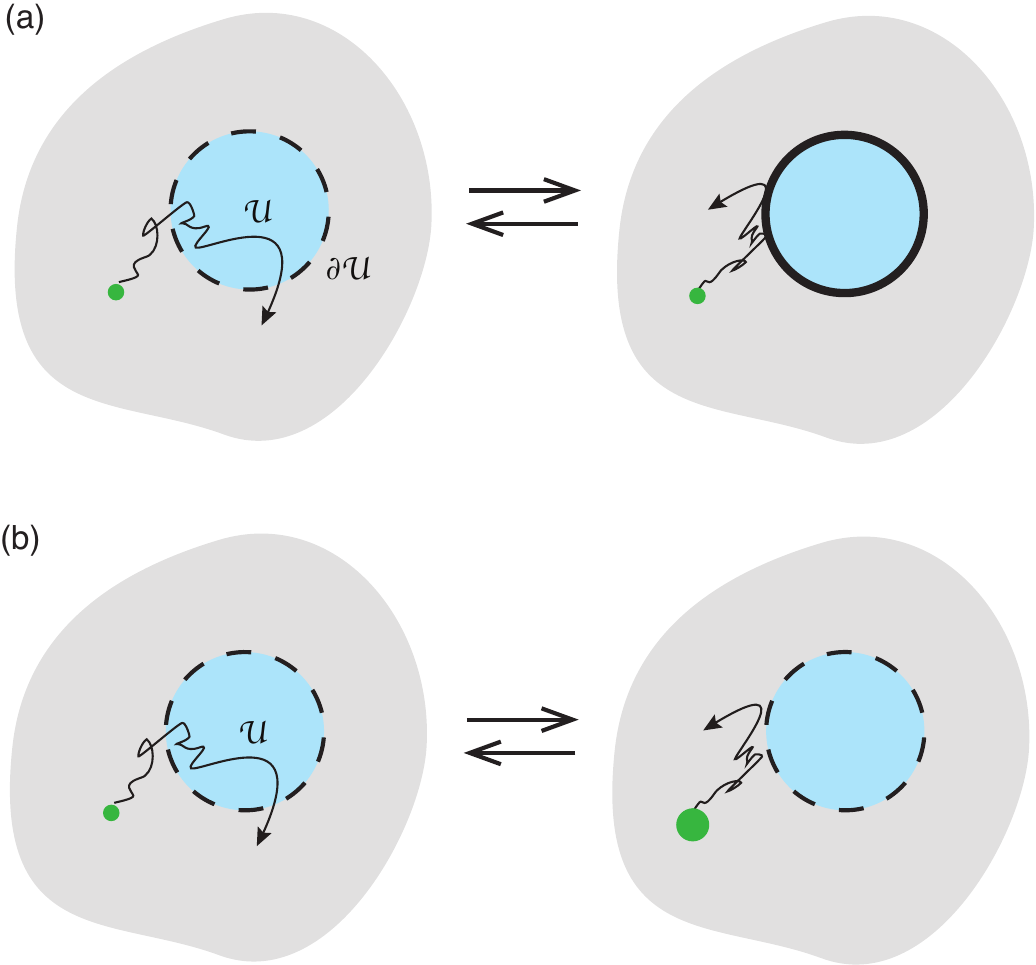}
  \caption{Particle versus boundary switching. (a) The boundary $\partial \calU$ switches between an open and a closed state. (b) The particle switches between two configuration states, only one of which can cross $\partial \calU$.}
  \label{fig9}
  \end{center}
\end{figure}

\noindent \underline{Switching boundaries versus switching particles.} One significant difference between the model of volume transmission and the model of active phase separation is that in the former case the boundaries physically switch, whereas in the latter case it is the particles that independently switch states. These two distinct scenarios are. illustrated in Fig. \ref{fig9}. For a sufficiently large number of solute molecules, the diffusion equations of active phase separation are deterministic. On the other hand, all particles move in the same randomly switching environment during volume transmission, which means that equations (\ref{VT}) are stochastic. In other words $u(\x,t)$ is a stochastic field, and equations (\ref{VT2}) are first-moment equations. It is also possible to derive singularly-perturbed diffusion equations for the second-order moments
$C_n(\x_1,\x_2,t)=\E[u(\x_1,t)u(\x_2,t){\bm 1}_{N(t)=n}]$
etc. From these it can be shown that the randomly switching environment induces statistical correlations between the particles even though they are not interacting. For example, 
\begin{equation}
\E[u(\x_1,t)u(\x_2,t){\bm 1}_{N(t)=n}]\neq  \E[u(\x_1,t){\bm 1}_{N(t)=n}]\E[u(\x_2,t){\bm 1}_{N(t)=n}].
\end{equation}
One direction for future exploration would be to apply singular perturbation theory to these higher-order moment equations.

\subsection{The Kuramoto model} The original non-spatial model (\ref{mass0}) has motivated a number of theoretical studies of simplified oscillator models coupled via quorum sensing \cite{Schwab12,Schwab12a,Sharma16,Verma19}. In particular, following Ref. \cite{Schwab12,Schwab12a}, suppose that each oscillator is just beyond a Hopf bifurcation point so that its state at time $t$ can be represented by a complex amplitude $z_j(t)$, $j=1,\ldots,N$\footnote{This contrasts with the regime considered in Ref. \cite{Iyan21b}, where each isolated oscillator is below the Hopf bifurcation threshold.}. It follows that the state of the environment can be represented by a complex variable $Z(t)$, which diffusively couples to each $z_j$. Let $\omega_j$ be the natural frequency of the $j$-th oscillator. The frequencies $\omega_j$ are randomly drawn from a distribution $h(\omega)$, which is taken to be an even function about a mean frequency $\omega_0$.
Finally, moving to a rotating frame with frequency $\omega_0$, the dynamics of the system can be represented by the equations 
\begin{subequations}
\label{QSz}
\begin{align}
\frac{dz_j}{dt}&=(\lambda_0+i\omega_j-|z_j|^2)z_j-\widehat{\kappa}(z_j-Z),\\
\frac{dZ}{dt} &=\frac{\alpha \widehat{\kappa}}{N}\sum_{j=1}^N (z_j-Z)-(\gamma_0 +i\omega_0)Z.
\end{align}
\end{subequations}
In the rotating frame the frequencies $\omega_j$ are drawn from an even distribution $g(\omega)=h(\omega-\omega_0)$ with zero mean. Suppose that equation (\ref{QSz}a) is rewritten in polar coordinates with $z_j=a_j\e^{i\theta_j}$ and $z=a\e^{i\psi}$:
\begin{subequations}
\begin{align}
\frac{da_j}{dt}&=(\lambda_0-\widehat{\kappa}-a_j^2)a_j +\widehat{\kappa}a \cos (\psi-\theta_j),\\
\frac{d\theta_j}{dt}&=\omega_j +\frac{\widehat{\kappa}a}{a_j}\sin(\psi-\theta_j).
\end{align}
\end{subequations}
As a further simplification, consider the large $\lambda_0$ limit with $a_j\approx \sqrt{\lambda_0}$ in steady-state for all $j=1,\ldots,N$. This allows us to effectively eliminate the dynamics of the amplitudes $a_j$. After performing the rescalings $a_j \rightarrow a_j/\sqrt{\lambda_0}$ and $a\rightarrow a/\sqrt{\lambda_0}$, equations (\ref{QSz}) reduce to \cite{Schwab12a}
\begin{subequations}
\label{Kur}
\begin{align}
\frac{d\theta_j}{dt}&=\omega_j+\widehat{\kappa}a\sin (\psi-\theta_j),\\
\frac{dz}{dt}&=\frac{\alpha \widehat{\kappa}}{N}\sum_{j=1}^N  (\e^{i\theta_j}-z)-(\gamma_0+i\omega_0)z.
\end{align}
\end{subequations}
Equations (\ref{Kur}) are a modified version of the classical Kuramoto model \cite{Kuramoto84,Strogatz00}, in which $z$ is identified with the state of an external medium rather than the first circular moment of the $N$ oscillators. The latter is defined according to
\begin{equation}
\bar{z}(t):=\frac{1}{N}\sum_{j=1}^N\e^{i\theta_j(t)}=\overline{a}(t)\e^{i\overline{\psi}(t)}.
\end{equation}
Here $\overline{\phi}(t)$ is equal to the average phase and $\overline{a}(t)$ is a measure of the degree of
phase-coherence;
a completely incoherent state corresponds to the case $\overline{a}=0$, whereas a completely synchronised state satisfies $\overline{a}=1$. 

The collective behaviour of the globally-coupled model given by equations (\ref{Kur}) can be investigated in the thermodynamic limit $N\rightarrow \infty$ by considering a continuum mean field model \cite{Schwab12a}. Let $\rho(\theta,t,\omega)d\theta$ denote the density of oscillators with frequency $\omega$ that have a phase in the interval $[\theta,\theta+d\theta]$ at time $t$ with the normalization
 \begin{equation}
\int_0^{2\pi}\rho(\theta,t,\omega)d\theta =g(\omega). 
\end{equation} 
Since the total number of oscillators with frequency $\omega$ is fixed, we have the continuity or Liouville equation
\begin{subequations}
\label{cont}
\begin{equation}
\frac{\partial u}{\partial t}=-\frac{\partial}{\partial \theta}\left [\omega u+\frac{\widehat{\kappa}}{2 i}(z\e^{-i\theta}-z^* \e^{i\theta})u \right ],
\end{equation}
and
\begin{equation}
\frac{dz}{dt}= {\alpha \widehat{\kappa}}  \left (\overline{z} -z\right )-(\gamma+i\omega_0)z,
\end{equation}
\end{subequations}
where
\begin{equation}
\overline{z}=\langle \e^{i\theta}\rangle :=\int_{-\infty}^{\infty}  \left [ \int_0^{2\pi}\e^{i\theta}u(\theta,t,\omega)\frac{d\theta}{2\pi}\right ]d\omega.
\end{equation}

The modified ODE system (\ref{wu3}), which takes into account the effects of finite diffusion in the bulk, suggests a generalisation of equations (\ref{Kur}) of the following form:
\begin{subequations}
\label{Kur2}
\begin{align}
\frac{d\theta_j}{dt}&=\omega_j+\widehat{\kappa}a\sum_{j=1}^N W_{jk}\sin (\psi-\theta_k),\\
\frac{dz}{dt}&=\frac{\alpha \widehat{\kappa}}{N}\sum_{j,k=1}^N W_{jk} (\e^{i\theta_k}-z)-(\gamma_0+i\omega_0)z.
\end{align}
\end{subequations}
It would be interesting to explore how the presence of the matrix ${\bf W}$ affects synchronization and clustering of the original model for finite $N$. However, in contrast to equations (\ref{Kur}), taking the thermodynamic limit of equations (\ref{Kur2}) doesn't make sense, since the reduction of the PDE-ODE system to the ODE system (\ref{wu3}) assumes that the compartments are well-separated. This assumption would break down in the limit $N\rightarrow \infty$, at least in a bounded domain $\Omega$.

\section*{Acknowledgements} I would like to thank Michael Ward (University of British Columbia) for many stimulating discussions regarding the asymptotic methods that he and his collaborators pioneered. 
	   
    \setcounter{equation}{0}
\renewcommand{\theequation}{A.\arabic{equation}}
 
\section*{Appendix A: Neumann Green's functions in simple geometries}

One of the crucial features of the asymptotic expansions presented in this review is the involvement of the Neumann Green's function $G(\x,\x')$ of the modified Helmholtz equation and the generalised Neumann Green's function of Laplace's equation, see equations (\ref{GMH}) and (\ref{exG1}), respectively. Explicit expressions for these Green's functions and their regular parts are only known for a small set of geometries, particularly in the case of the modified Helmholtz equation. In this appendix we collect some of the known results.
\medskip

\noindent (a) {\em Laplace's equation on the disc}. Let $\calU \subset \R^2$ be the unit circle centred at the origin. The 2D Neumann Green's function is given by
\begin{align}
G_0(\x,\hxi)&=\frac{1}{2\pi}\left [-\ln(|\x-\hxi|)-\ln\left (\left |\x|\hxi|-\frac{\hxi}{|\hxi |}\right |\right )  +\frac{1}{2}(|\x|^2+|\hxi|^2)-\frac{3}{4}\right ],
\label{Gdisc}
\end{align}
with the regular part obtained by dropping the first logarithmic term. 
\medskip

\noindent (b) {\em Laplace's equation the sphere.} Let $\calU \subset \R^3$ be the sphere of radius $a$ centred about the origin. The 3D Neumann Green's function takes the form \cite{Cheviakov11}
\begin{align}
G_0(\x,\hxi)&=\frac{1}{4\pi |\x-\hxi|}+\frac{a}{4\pi |\x|r'} +\frac{1}{4\pi a}\ln \left (\frac{2a^2}{a^2-|\x||\hxi|\cos \theta+|\x|r'}\right ) \nonumber \\
&\quad 
 +\frac{1}{6|\calU|}(|\x|^2+|\hxi|^2)+B,
\end{align}
where the constant $B$ is chosen so that $\int_{\calU}G(\x,\hxi)d\x=0$, and
\[\cos \theta =\frac{\x\cdot \hxi}{|\x||\hxi|},\quad \x'=\frac{a^2\x}{|\x|^2},\quad r'=|\x'-\hxi|.
\]
It can be shown that $B$ is independent of $\hxi$.
\medskip

\noindent (c) {\em Laplace's equation in a rectangular domain.} Let $\calU \subset \R^2$ be a rectangular domain $[0,L_1]\times [0,L_2]$. The 2D Neumann Green's function has the logarithmic expansion \cite{Bressloff08,Pillay10}
\begin{align}
  G_0 (\r,\r') &=\frac{1}{L_{1}}H_{0}(y,y')  -\frac{1}{2\pi } \sum_{j=0}^{\infty}\sum_{n = \pm}\sum_{m=\pm} \left(\ln\abs{1-\tau^{j}z_{n}\zeta_{m}}+\ln\abs{1-\tau^{j}z_{n}\varsigma_{m}}\right),
\label{eq6:G0}
\end{align}
where
\begin{equation*}
\tau= e^{-2\pi L_{2}/L_{1}},\quad z_{\pm} = \e^{i\pi(x\pm x')/L_{1}},
\end{equation*}
\begin{equation*}
 \zeta_{\pm}=e^{-\pi\abs{y\pm y'}/L_{1}} ,\quad \varsigma_{\pm} = e^{-\pi(2L_{2}-\abs{y\pm y'})/L_{1}},
\end{equation*}
and
\begin{equation}
  \label{H0}
  H_{0}(y,y') = \frac{L_{2}}{3} + \frac{1}{2L_{2}}(y^{2}+{y'}^{2}) - \max\{y,y'\}.
\end{equation}
Assuming that $\tau \ll 1$, we have the approximation
\begin{align*}
G(\r,\r') &=\frac{1}{L_{1}}H_{0}(y,y') -\frac{1}{2\pi } \sum_{n = \pm}\sum_{m=\pm} \left(\ln\abs{1-z_{n}\zeta_{m}}+\ln\abs{1-z_{n}\varsigma_{m}}\right) +O(\tau) .
\end{align*}
The only singularity exhibited by equation (\ref{eq6:G0}) occurs when $\r\to \r'$, $\r'\notin \partial \calU$, in which case $z_-=\zeta_-=1$ and the term $\ln\abs{1-z_{-}\zeta_{-}}$ diverges. Writing
\begin{equation*}
   \ln\abs{1-z_{-}\zeta_{-}} = \ln\abs{\r-\r'} + \ln\frac{\abs{1-z_{-}\zeta_{-}}}{\abs{\r-\r'}},
\end{equation*}
where the first term on the right hand side is singular and the second is regular, we find that
\begin{equation}
  G (\r,\r') = -\frac{1}{2\pi}\ln\abs{\r-\r'} +R(\r,\r'),
\end{equation}
where $R$ is the regular part of the Green's function given by 
\begin{align}
  R(\r,\r') &=-\frac{1}{L_{1}} H_{0}(y,y')+\frac{1}{2\pi}\ln\frac{\abs{1-z_{-}\zeta_{-}}\abs{1-z_{-}\zeta_{+}}}{\abs{\r-\r'}}\nonumber \\
&\quad +\frac{1}{2\pi}\ln\abs{1-z_{-}\varsigma_{-}}\abs{1-z_{-}\varsigma_{+}}+ \frac{1}{2\pi}\ln\abs{1-z_{+}\varsigma_{-}}\abs{1-z_{+}\varsigma_{+}}\nonumber \\
&\quad +\frac{1}{2\pi}\ln\abs{1-z_{+}\zeta_{-}}\abs{1-z_{+}\zeta_{+}} + O(\tau) .
\end{align}

\noindent {\em (d) Modified Helmholtz equation in a sphere \cite{Grebenkov20}.}
\begin{align}
{G} (\x,s|\x_0)&=\frac{\e^{-\sqrt{s/D}|\x-\x_0|}}{4\pi |\x-\x_0|}-{G}_{\rm sp}(\x,s|\x_0),
\end{align}
with
\begin{align}
{G}_{\rm sp}(\x,s|\x_0)&=\frac{1}{4\pi}\sqrt{\frac{s}{D}}\sum_{n=0}^{\infty} (2n+1)P_n(\cos \theta)\frac{k_n'(\sqrt{s/D}R_0)}{i_n'(\sqrt{s/D}R_0)}\nonumber \\
&\qquad \times i_n(\sqrt{s/D}|\x|)
i_n(\sqrt{s/D}|\x_0|).
\label{Gsphere}
\end{align}
Here $P_n$ is a Legendre polynomial, $\x\cdot \x_0 =|\x||\x_0|\cos \theta$, and $i_n,k_n$ are modified spherical Bessel functions,
\begin{equation}
i_n(x)=\sqrt{\frac{\pi}{2x}}I_{n+1/2}(x),\quad k_n(x)=\sqrt{\frac{2}{\pi x}}K_{n+1/2}(x).
\end{equation}
In order to simplify the analysis, we will assume that the initial position $\x_0$ is at the centre of the sphere. Using the identities
\begin{equation}
i_0(x)=\frac{\sinh x}{x},\quad k_0(x)=\frac{\e^{-x}}{x}, \quad i_n(0)=0, \ n>0,
\end{equation}
we see that
\begin{align}
G_{j0}&=\frac{e^{-\sqrt{s/D}|\x_j|}}{4\pi |\x_j|}-\frac{1}{4\pi}\sqrt{\frac{s}{D}}\frac{k_0'(\sqrt{s/D}R_0)}{i_0'(\sqrt{s/D}R_0)}i_0(\sqrt{s/D}|\x_j|).
\label{Gsphere2}
\end{align}
Similarly,
\begin{equation}
{\mathcal G}_{jk}=\frac{e^{-\sqrt{s/D}|\x_j-\x_k|}}{4\pi |\x_j-\x_k|}-{G}_{\rm sp}(\x_j,s|\x_k),\ j\neq k; \quad {\mathcal G}_{jj}=-\frac{1}{4\pi}\sqrt{\frac{s}{D}}-G_{\rm sp}(\x_j,s|\x_j).
\end{equation}
Further useful identities are
\begin{equation}
(2n+1) i'_n=ni_{n-1}+(n+1)i_{n+1},\quad -(2n+1) k'_n=nk_{n-1}+(n+1)k_{n+1}.
\end{equation}
The 3D Neumann Green's function for Laplace's equation in the sphere can also be written down explicitly so that \cite{Cheviakov11} 
\begin{align}
\overline{G}_{jk}=\overline{G}(\x_j,\x_k)&=\frac{1}{4\pi }\bigg [\frac{1}{|\x_j-\x_k|}+\frac{R_0}{|\x_j||\x_j'-\x_k|} +\frac{1}{2}(|\x_j|^2+|\x_k|^2) \\
&\quad +\frac{1}{R_0}\ln \left (\frac{2R_0^2}{R_0^2-|\x_j||\x_k|\cos \theta+|\x_j||\x_j'-\x_k|}\right )\bigg ]-\frac{7}{10 \pi R_0} \nonumber 
\end{align}
for $j\neq k$, where $\x'=\x/|\x|^2$.
It also follows that
\begin{align}
\overline{G}_{jj}=\overline{R}(\x_j,\x_j)&=\frac{1}{4\pi }\bigg [\frac{R_0}{R_0^2-|\x_j|^2}+\frac{1}{R_0}\ln \left (\frac{R_0^2}{R_0^2-|\x_j|^2}\right )+|\x_j|^2 \bigg ]-\frac{7}{10 \pi R_0},
\end{align}
and
\begin{equation}
\overline{G}_{j0}=\frac{1}{4\pi}\left [\frac{1}{|\x_j|}+\frac{|\x_j|^2}{2}\right ] +1-\frac{7}{10 \pi R_0}.
\end{equation}

\end{document}